\documentclass[11pt]{amsart}
\usepackage[T1]{fontenc}
\usepackage[utf8]{inputenc}
\usepackage{lmodern}
\usepackage{amsmath,amssymb,mathtools}
\usepackage[a4paper,margin=28mm]{geometry}
\usepackage{microtype}
\usepackage[shortlabels]{enumitem}
\usepackage{booktabs}
\usepackage{array}
\usepackage{tikz}
\usepackage[hidelinks]{hyperref}

\newtheorem{theorem}{Theorem}[section]
\newtheorem{proposition}[theorem]{Proposition}
\newtheorem{lemma}[theorem]{Lemma}
\newtheorem{corollary}[theorem]{Corollary}

\newtheorem{introthm}{Theorem}

\newtheorem{conjecture}[theorem]{Conjecture}
\theoremstyle{definition}
\newtheorem{definition}[theorem]{Definition}
\newtheorem{example}[theorem]{Example}

\theoremstyle{remark}
\newtheorem{remark}[theorem]{Remark}
\numberwithin{equation}{section}

\newcommand{\R}{\mathbb R}
\newcommand{\st}{\mathord{*}}
\newcommand{\cF}{\mathcal F}
\newcommand{\cP}{\mathcal P}
\DeclareMathOperator{\supp}{supp}
\DeclareMathOperator{\aff}{aff}
\DeclareMathOperator{\relint}{relint}
\DeclareMathOperator{\interior}{int}
\DeclareMathOperator{\conv}{conv}

\title[More than $2^d$ neighbourly simplices]{More than $2^d$ neighbourly
simplices in $\mathbb R^d$}
\author[J. Grytczuk]{Jaros\l aw Grytczuk}
\address{Faculty of Mathematics and Information Science,
Warsaw University of Technology, 00-662 Warsaw, Poland}
\email{jaroslaw.grytczuk@pw.edu.pl}
\author[A. P. Kisielewicz]{Andrzej P. Kisielewicz}
\address{Instytut Matematyki, Informatyki i Ekonometrii,
Uniwersytet Zielonog\'orski, ul. Podg\'orna 50,
65-246 Zielona G\'ora, Poland}
\email{A.Kisielewicz@wmie.uz.zgora.pl}
\author[K. Przes\l awski]{Krzysztof Przes\l awski}
\address{Instytut Matematyki, Informatyki i Ekonometrii,
Uniwersytet Zielonog\'orski, ul. Podg\'orna 50,
65-246 Zielona G\'ora, Poland}
\email{K.Przeslawski@wmie.uz.zgora.pl}
\date{}

\begin{document}
\begin{abstract}
Two $d$-simplices in $\R^d$ are neighbourly if their intersection is
$(d-1)$-dimensional. Let $S_d$ be the largest size of a family of pairwise
neighbourly $d$-simplices in $\R^d$. Zaks constructed $2^d$ such simplices,
and it has been conjectured that $S_d=2^d$; this is true for $d\le3$. We
construct $18$ pairwise neighbourly $4$-simplices in $\R^4$. Combined with
Zaks's doubling step this gives $S_d\ge 9\cdot2^{d-3}$ for every $d\ge4$, so
the conjecture fails in every dimension $d\ge4$. The construction rests on
nine $4$-simplices that lie on one side of a hyperplane and have a facet in
it. Nine is one more than the largest number of pairwise neighbourly
tetrahedra, and indeed eleven of the $36$ contacts take place strictly above
the hyperplane.

For upper bounds we use neighbourly codes: sets of words of weight $k$ over
$\{0,1,\st\}$ with pairwise different supports, every two of which conflict
in exactly one position. Their largest size $M_k$ satisfies
$S_d\le M_{d+1}$. We prove $M_k<\tfrac{35}{46}\,2^k$ for all $k$ and
$M_k=3\cdot2^{k-2}$ for $k\le8$. Consequently $S_d<\tfrac{35}{23}\,2^d$ for
all $d$, and $18\le S_4\le24$. We also show that a neighbourly family in
$\R^d$ with more than $2^d$ members contains two rainbow triples sharing two
members.
\end{abstract}
\maketitle

\section{Introduction}\label{sec:intro}
Two $d$-dimensional simplices in $\R^d$ are \emph{neighbourly} if their
intersection is $(d-1)$-dimensional, and a family of $d$-simplices is
\emph{neighbourly} if every two of its members are neighbourly. Let $S_d$
denote the largest cardinality of a neighbourly family of $d$-simplices in
$\R^d$.

The problem goes back to Bagemihl~\cite{Bag}, who proved $8\le S_3\le17$ and
conjectured $S_3=8$. Baston~\cite{Bas} improved the upper bound to $S_3\le9$,
and the conjecture $S_3=8$ was finally confirmed by Zaks~\cite{Z91}. Earlier
Zaks~\cite{Z81} had constructed $2^d$ pairwise neighbourly $d$-simplices in
$\R^d$ for every $d$. The general conjecture
\[
  S_d=2^d\qquad(d\ge1)
\]
is attributed to Bagemihl and Baston (see~\cite{Kal,K24}), and it is known
to hold for $d\le3$. Perles~\cite{Per}
found a clever proof of $S_d\le2^{d+1}$ (see also~\cite{AZ,Kal});
this was improved to $2^{d+1}-1$ in~\cite{AZ} and to $2^{d+1}-2$ by the second
and third authors~\cite{KP19,KP21}, and Kisielewicz~\cite{K24} proved that
$2^{d+1}-S_d\to\infty$. Neighbourly families of boxes, where the
combinatorics is much better understood, were studied by Alon~\cite{Alon} and
in~\cite{AGKP}.

\subsection*{Lower bound}
Our first result shows that the conjecture fails in every dimension $d\ge4$.

\begin{introthm}\label{thm:A}
There exist $18$ pairwise neighbourly $4$-simplices in $\R^4$. Consequently
\[
   S_d\ \ge\ 18\cdot2^{d-4}\;=\;\tfrac98\,2^d\qquad\text{for every }d\ge4 .
\]
\end{introthm}

Zaks's construction can be described as follows. Place a neighbourly family
of $(d-1)$-simplices in a hyperplane $H\subset\R^d$ and erect pyramids over
them with a common apex on one side of $H$. Take a suitably transformed copy
on the other side of $H$, arranged so that every base from one side overlaps
every base from the other side. Each half of such a family consists of
simplices lying on one side of $H$ and having a facet in $H$; we call such
families \emph{based}. The half produced by pyramids has at most $S_{d-1}$
members.

The family of Theorem~\ref{thm:A} is obtained in the same way from a based
family of \emph{nine} pairwise neighbourly $4$-simplices
(Section~\ref{sec:nine}). Nine exceeds $S_3=8$, so these simplices are not
pyramids over a neighbourly family of tetrahedra: for $11$ of the $36$ pairs
the bases are disjoint, and the two simplices touch only above $H$
(Remark~\ref{rem:heights}). The doubling step itself, due to
Zaks~\cite{Z81}, is proved in Section~\ref{sec:doubling} in the form we need.
All numerical claims about the construction can be checked by hand from the
certificates printed in Section~\ref{sec:nine}. They have also been verified
in exact rational arithmetic by two independent programs, one of which
accompanies this paper.

\subsection*{Upper bounds}
We use the coding of neighbourly families introduced
in~\cite{KP19,K24}. A \emph{neighbourly $k$-code} is a finite set of words
over $\{0,1,\st\}$, each with exactly $k$ letters different from $\st$, such
that any two words have different supports and differ in exactly one position
where both have a letter from $\{0,1\}$ (precise definitions are in
Section~\ref{sec:codes}). Let $M_k$ be the largest size of a neighbourly
$k$-code. The code of a neighbourly family of $d$-simplices is a neighbourly
$(d+1)$-code, whence $S_d\le M_{d+1}$ (Lemma~\ref{lem:coding}). A simple
doubling gives $M_k\ge g_k:=3\cdot2^{k-2}$. The second and third authors
showed $M_k\le2^k-2$~\cite{KP19} and conjectured that $M_k=g_k$
\cite[Sect.~3]{KP19}; Kisielewicz~\cite{K24} proved $2^k-M_k\to\infty$.

\begin{introthm}\label{thm:B}
For every $k\ge2$,
\[
  g_k\ \le\ M_k\ \le\ \max\Bigl\{g_k,\ \Bigl\lfloor
     \frac{35\cdot2^k}{46+78/g_k}\Bigr\rfloor\Bigr\},
  \qquad g_k=3\cdot2^{k-2}.
\]
In particular $M_k<\frac{35}{46}\,2^k$ for every $k$, and $M_k=g_k$ for
$2\le k\le8$.
\end{introthm}

\begin{introthm}\label{thm:C}
For every $d\ge1$ we have $S_d<\frac{35}{23}\,2^d$. Moreover,
\[
  18\le S_4\le24,\qquad 36\le S_5\le48,\qquad 72\le S_6\le96,\qquad
  144\le S_7\le192 .
\]
\end{introthm}

The equality $M_5=24$, which gives $S_4\le24$, already follows from the first
part of the argument (Corollary~\ref{cor:M5}). If the extremal $5$-code is
unique up to isomorphism, then $S_4\le23$ (Proposition~\ref{prop:S4}).

\subsection*{Rainbow triples}
For two members $P,Q$ of a neighbourly family, $P\cap Q$ spans a hyperplane
$\aff(P\cap Q)$, the \emph{colour} of the pair. Three members form a
\emph{rainbow triple} if the three colours of their pairs are distinct.

\begin{introthm}\label{thm:D}
Let $\cF$ be a neighbourly family of $d$-simplices in $\R^d$. If no two
rainbow triples of $\cF$ have two common members, then $|\cF|\le2^d$. More
generally, if $|\cF|=2^d+u$ with $u>0$, then $\cF$ contains at least
$\lceil u/4\rceil$ pairwise disjoint quadruples $\{A,B,C,D\}$ such that
$\{A,B,C\}$ and $\{A,B,D\}$ are rainbow triples.
\end{introthm}

For the $18$ simplices of Theorem~\ref{thm:A}, $46$ of the $153$ pairs lie in
two or more rainbow triples (Example~\ref{ex:rainbow18}). The proof of
Theorem~\ref{thm:D} combines a geometric fact with a purely combinatorial
statement. The geometric fact is that the code of a neighbourly family in
$\R^d$ contains no three words sharing $d-1$ letters (Lemma~\ref{lem:line}).
The combinatorial statement is Theorem~\ref{thm:linearnocore}, a structural
result on codes whose rainbow triples form a linear hypergraph.

\subsection*{Organisation}
Section~\ref{sec:codes} collects definitions and basic facts about codes.
Section~\ref{sec:nine} presents the nine based simplices in $\R^4$, and
Section~\ref{sec:doubling} proves Theorem~\ref{thm:A}. Section~\ref{sec:rainbow}
studies rainbow triples and proves Theorem~\ref{thm:D}. Section~\ref{sec:upper}
proves Theorems~\ref{thm:B} and~\ref{thm:C}. Open problems are collected in
Section~\ref{sec:open}, and the vertices of the nine simplices are listed in
Appendix~\ref{app:vertices}.

\section{Neighbourly codes}\label{sec:codes}
\subsection{Words, conflicts and cylinders}
We consider words $v=v_1\cdots v_n$ over the alphabet $\{0,1,\st\}$. The
\emph{support} of $v$ is $\supp v=\{i:v_i\ne\st\}$ and its \emph{weight} is
$|\supp v|$; the letters $0,1$ are called \emph{bits}. Two words $u,v$
\emph{conflict} at position $i$ if $\{u_i,v_i\}=\{0,1\}$.

\begin{definition}\label{def:code}
A finite set $V$ of words of length $n$ is a \emph{neighbourly $k$-code} if
every word of $V$ has weight $k$, every two distinct words of $V$ conflict at
exactly one position, and distinct words of $V$ have distinct supports.
\end{definition}

In the terminology of~\cite{K24} these are the neighbourly $k$-codes without
twin pairs; two words with equal supports and exactly one conflict form a
\emph{twin pair}. We identify codes that differ by a permutation of
positions, by exchanging $0$ and $1$ in some columns, or by adding or
deleting columns consisting of stars only; this is \emph{isomorphism}. Let
$M_k$ be the largest size of a neighbourly $k$-code (of any length).

For $u\ne v$ in $V$ the unique conflict position $c(u,v)$ is the
\emph{colour} of the pair $uv$. For a position $i$ and $a\in\{0,1\}$ let
$V^{i,a}=\{u\in V:u_i=a\}$. The column $i$ is \emph{mixed} if both
$V^{i,0}$ and $V^{i,1}$ are nonempty; the unordered pair
$\{|V^{i,0}|,|V^{i,1}|\}$ is its \emph{profile}, written $(a,b)$ with $a\ge b$.

Let $\mu$ be the uniform probability measure on $\{0,1\}^n$. The
\emph{cylinder} of $v$ is
$[v]=\{x\in\{0,1\}^n:x_i=v_i\text{ for all }i\in\supp v\}$; it has measure
$2^{-|\supp v|}$. Two words conflict if and only if their cylinders are
disjoint. Hence the cylinders of a neighbourly $k$-code are pairwise
disjoint, $M_k\le2^k$, and the maximum $M_k$ is attained. For a code $V$ of
weight $k$ we measure volume in \emph{units} of $2^{-k}$, so that the cube has
volume $2^k$ and every cylinder of a word of $V$ has volume $1$. The number
\[
   \Delta(V)=2^k-|V|
\]
is the volume not covered by the cylinders of $V$.

\subsection{Doubling and the codes \texorpdfstring{$G_k$}{G\_k}}
For a code $G$ of length $s$ let
\[
  D(G)=\{0\,u\,\st^s:u\in G\}\cup\{1\,\st^s\,u:u\in G\}.
\]
Words on the same side keep their conflicts, words from different sides
conflict only at the first position, and all supports are distinct.
Therefore $D(G)$ is a neighbourly $(k+1)$-code if $G$ is a neighbourly
$k$-code, and
\begin{equation}\label{eq:Mdouble}
   M_{k+1}\ge2M_k .
\end{equation}
Let
\[
   T=\{00\st,\ 1\st0,\ \st11\},
\]
a neighbourly $2$-code, and put $G_k=D^{k-2}(T)$ for $k\ge2$. Then $G_k$ is a
neighbourly $k$-code with $g_k=3\cdot2^{k-2}$ words, and the cylinders of
$G_k$ cover three quarters of the cube.

\begin{conjecture}[\cite{KP19}]\label{conj:34}
$M_k=g_k$ for every $k\ge2$.
\end{conjecture}

By~\eqref{eq:Mdouble} the sequence $M_k/2^k$ is nondecreasing; hence
Conjecture~\ref{conj:34} is equivalent to $\lim_k M_k/2^k=3/4$.

\subsection{From simplices to codes}
Let $\cF$ be a finite neighbourly family of $d$-simplices in $\R^d$, $d\ge1$.
Let $H_1,\dots,H_n$ be the distinct hyperplanes spanned by facets of members
of $\cF$, and for each $i$ label the two closed half-spaces bounded by $H_i$
by $0$ and $1$. The word $w(P)\in\{0,1,\st\}^n$ of $P\in\cF$ has
$w(P)_i=a$ if $H_i$ spans a facet of $P$ and $P$ lies in the half-space
labelled $a$, and $w(P)_i=\st$ otherwise.

\begin{lemma}[\cite{KP19,K24}]\label{lem:coding}
The map $w$ is injective and $w(\cF)$ is a neighbourly $(d+1)$-code. For
$P\ne Q$ in $\cF$, the colour of $w(P)w(Q)$ is the index of the hyperplane
$\aff(P\cap Q)$. In particular $S_d\le M_{d+1}$.
\end{lemma}

\begin{proof}
Every word has weight $d+1$. Let $P\ne Q$ be members of $\cF$. Since
$P\cap Q$ is $(d-1)$-dimensional, the interiors of $P$ and $Q$ are disjoint,
and $P$ and $Q$ are separated by a hyperplane $H$, which must contain
$P\cap Q$; thus $H=\aff(P\cap Q)$. The face $P\cap H$ of $P$ contains the
$(d-1)$-dimensional set $P\cap Q$, so it is a facet of $P$; likewise for $Q$.
Hence $H=H_i$ for some $i$ and $w(P),w(Q)$ conflict at $i$. A second conflict
at $j\ne i$ would give $P\cap Q\subseteq H_i\cap H_j$, which has dimension
at most $d-2$. Finally, if $\supp w(P)=\supp w(Q)$, then $Q$ is obtained from
the system of facet inequalities of $P$ by reversing exactly one inequality,
say the one of the facet opposite to the vertex $p$ of $P$. The remaining $d$
inequalities define a pointed cone with apex $p$, and reversing the last one
leaves the unbounded part of this cone beyond the facet, which is not a
simplex. This contradiction proves that the supports are distinct; in
particular $w$ is injective.
\end{proof}

\subsection{Triangles and rainbow triples}
\begin{definition}\label{def:triangle}
Let $V$ be a neighbourly $k$-code, $k\ge2$. A \emph{triangle} in $V$ is a set
of three words of $V$ having equal bits on a set of $k-2$ positions. The
pattern formed by these common bits (with stars elsewhere) is the \emph{core}
of the triangle, and the remaining positions of the supports are its
\emph{local} positions. A set of three words is a \emph{rainbow triple} if
the three colours of its pairs are distinct.
\end{definition}

After the core is deleted, the three words of a triangle have weight two and
distinct supports, and every two of them conflict. Their two-element
supports pairwise intersect; they cannot share one position, because three
words cannot pairwise conflict at the same position. So the supports are the
three edges of a triangle on three local positions $p,q,r$, and up to
isomorphism the words are $00\st,\ 1\st0,\ \st11$ on $(p,q,r)$, that is, a
copy of $T$. In particular every triangle is a rainbow triple whose colours
are its local positions.

In a rainbow triple, at the conflict position of any two of its words the
third word has a star (otherwise two colours would coincide), and at
positions in the common support of all three words their bits agree.
Consequently a rainbow triple is a triangle if and only if the common
support of its three words has $k-2$ elements.

\section{Nine simplices on one side of a hyperplane}\label{sec:nine}
We use coordinates $(x,y,z,t)$ in $\R^4$. For $1\le k\le11$ let
\[
   F_k(x,y,z,t)=x+p_k\,y+q_k\,z+v_k\,t-h_k ,
\]
with the integer coefficients given in Table~\ref{tab:forms}. For
$1\le i\le9$ the set $P_i$ is defined by $t\ge0$ and four inequalities
$\varepsilon F_k\ge0$, $\varepsilon\in\{+1,-1\}$, listed in row $i$ of
Table~\ref{tab:words}; there $+k$ stands for $F_k\ge0$ and $-k$ for
$F_k\le0$. Equivalently, row $i$ is the word $w_i\in\{0,1,\st\}^{11}$ with
$1$ for $+k$, $0$ for $-k$ and $\st$ for the forms not used by $P_i$. The
last column gives a point $(u_i,0,0,0)$.

\begin{table}[ht]
\centering\small
\caption{The forms $F_k=x+p_ky+q_kz+v_kt-h_k$.}\label{tab:forms}
\begin{tabular}{@{}c*{11}{r}@{}}
\toprule
$k$ & 1 & 2 & 3 & 4 & 5 & 6 & 7 & 8 & 9 & 10 & 11\\
\midrule
$p_k$ & 14 & 5 & $-212$ & 40 & $-89$ & $-47$ & $-27$ & $-71$ & $-52$ & 26 & $-187$\\
$q_k$ & 143 & 171 & 236 & $-54$ & $-31$ & $-21$ & 10 & $-44$ & $-31$ & $-29$ & $-77$\\
$v_k$ & 0 & $-103$ & 0 & $-38$ & $-13$ & 2 & 0 & $-71$ & $-76$ & $-95$ & $-41$\\
$h_k$ & 38 & 4 & 18 & 24 & 41 & 39 & 37 & $-1$ & 1 & 3 & 5\\
\bottomrule
\end{tabular}
\end{table}

\begin{table}[ht]
\centering\small
\caption{The side facets of $P_i$, the corresponding word $w_i$ (positions
$1,\dots,11$), and the point $(u_i,0,0,0)$.}\label{tab:words}
\begin{tabular}{@{}c>{$}l<{$}c*{11}{c}c@{}}
\toprule
$i$ & \text{side facets} & & \multicolumn{11}{c}{$w_i$} & $u_i$\\
\midrule
1 & +8,\,-9,\,-10,\,-11 && $\st$&$\st$&$\st$&$\st$&$\st$&$\st$&$\st$&1&0&0&0 & $0$\\
2 & -5,\,+6,\,+7,\,+11  && $\st$&$\st$&$\st$&$\st$&0&1&1&$\st$&$\st$&$\st$&1 & $40$\\
3 & +3,\,+4,\,-7,\,+11  && $\st$&$\st$&1&1&$\st$&$\st$&0&$\st$&$\st$&$\st$&1 & $61/2$\\
4 & -2,\,-4,\,-7,\,+10  && $\st$&0&$\st$&0&$\st$&$\st$&0&$\st$&$\st$&1&$\st$ & $7/2$\\
5 & +2,\,-3,\,-7,\,+11  && $\st$&1&0&$\st$&$\st$&$\st$&0&$\st$&$\st$&$\st$&1 & $23/2$\\
6 & +1,\,-6,\,+7,\,+9   && 1&$\st$&$\st$&$\st$&$\st$&0&1&$\st$&1&$\st$&$\st$ & $77/2$\\
7 & -1,\,+9,\,-10,\,-11 && 0&$\st$&$\st$&$\st$&$\st$&$\st$&$\st$&$\st$&1&0&0 & $2$\\
8 & -1,\,+2,\,+10,\,-11 && 0&1&$\st$&$\st$&$\st$&$\st$&$\st$&$\st$&$\st$&1&0 & $9/2$\\
9 & -1,\,-6,\,+7,\,+11  && 0&$\st$&$\st$&$\st$&$\st$&0&1&$\st$&$\st$&$\st$&1 & $75/2$\\
\bottomrule
\end{tabular}
\end{table}

One checks directly that $W=\{w_1,\dots,w_9\}$ is a neighbourly $4$-code:
every two rows have exactly one position with opposite bits, and the
supports are distinct.

\begin{proposition}\label{prop:simplices}
Each $P_i$ is a $4$-simplex. Its \emph{base} $\beta_i=P_i\cap\{t=0\}$ is a
facet, and $P_i$ lies in the half-space $t\ge0$. At the point $(u_i,0,0,0)$
each of the four signed side forms $\varepsilon F_k$ of $P_i$ takes a value
at least $1/2$; in particular $(u_i,0,0,0)$ lies in the relative interior of
$\beta_i$.
\end{proposition}

\begin{proof}
Let $\varepsilon_1F_{k_1},\dots,\varepsilon_4F_{k_4}$ be the signed side
forms of $P_i$ in the order of Table~\ref{tab:words}. With the positive
integers of Table~\ref{tab:lambda} one has the identity
\begin{equation}\label{eq:lambda}
   \sum_{j=1}^4\lambda_j\,\varepsilon_jF_{k_j}=a_i-b_i\,t
\end{equation}
of affine functions on $\R^4$. Moreover, the vectors $(1,p_{k_j},q_{k_j})$,
$j=1,2,3$, are linearly independent (a $3\times3$ determinant). Hence the
affine map
\[
  \Psi_i:\R^4\to\R^5,\qquad
  \Psi_i=\bigl(b_i t,\ \lambda_1\varepsilon_1F_{k_1},\ \dots,\
          \lambda_4\varepsilon_4F_{k_4}\bigr)
\]
has rank $4$. By~\eqref{eq:lambda} it maps $\R^4$ bijectively onto the
hyperplane $\{y\in\R^5:y_0+\dots+y_4=a_i\}$, and $P_i$ onto the simplex
$\{y\ge0,\ \sum y_j=a_i\}$, which is $4$-dimensional because $a_i>0$. Its
facets correspond to the five coordinate hyperplanes; one of them is the
base. The values at $(u_i,0,0,0)$ are computed directly from the tables.
\end{proof}

\begin{table}[ht]
\centering\small
\caption{Certificates for Proposition~\ref{prop:simplices}:
$\sum_j\lambda_j\varepsilon_jF_{k_j}=a_i-b_it$.}\label{tab:lambda}
\begin{tabular}{@{}crrrrrr@{}}
\toprule
$i$&$\lambda_1$&$\lambda_2$&$\lambda_3$&$\lambda_4$&$a_i$&$b_i$\\\midrule
1 & 3318 & 1461 & 881 & 976 & 12302 & 831 \\
2 & 1610 & 583 & 476 & 551 & 22906 & 495 \\
3 & 16069 & 52255 & 71626 & 3302 & 1090290 & 2121072 \\
4 & 779 & 9781 & 2275 & 12835 & 283530 & 767410 \\
5 & 52255 & 22976 & 66296 & 37017 & 2472415 & 6899962 \\
6 & 15 & 548 & 70 & 463 & 17749 & 36284 \\
7 & 79 & 886 & 487 & 320 & 5177 & 7951 \\
8 & 3634 & 3101 & 634 & 101 & 124291 & 375492 \\
9 & 3220 & 17713 & 19544 & 1389 & 83094 & 92375 \\
\bottomrule
\end{tabular}
\end{table}

\begin{proposition}\label{prop:nine}
The simplices $P_1,\dots,P_9$ are pairwise neighbourly.
\end{proposition}

\begin{proof}
Let $i<j$. Since $W$ is a neighbourly code, there is exactly one $k$ such
that $F_k$ occurs in both $P_i$ and $P_j$ with opposite signs. Then
$P_i\subseteq\{\pm F_k\ge0\}$ and $P_j\subseteq\{\mp F_k\ge0\}$, so
$P_i\cap P_j$ lies in the hyperplane $\{F_k=0\}$ and has dimension at most
$3$. It suffices to exhibit a point of $\{F_k=0\}$ with $t>0$ at which all
the remaining inequalities of $P_i$ and $P_j$ are strict: a neighbourhood of
this point in $\{F_k=0\}$ is then contained in $P_i\cap P_j$.

Such points are listed in Table~\ref{tab:witness}. From the integers
$D,Y,Z,T$ in the row of $(i,j)$ put
\[
   (y,z,t)=(Y,Z,T)/D,\qquad x=h_k-p_ky-q_kz-v_kt ,
\]
so that $F_k=0$. Since $DF_r=D(h_k-h_r)+(p_r-p_k)Y+(q_r-q_k)Z+(v_r-v_k)T$,
all the values $D\varepsilon F_r$ ($r\ne k$) of the signed side forms of
$P_i$ and $P_j$ at this point are integers; the column $m$ gives their
minimum. In every row $T>0$ and $m>0$.
\end{proof}

\begin{table}[ht]
\centering\small
\setlength{\tabcolsep}{3.4pt}
\caption{Contact points for the $36$ pairs $(i,j)$; $k$ is the separating
form and $m$ is the minimum of $D\varepsilon F_r$ over the other signed
side forms.}\label{tab:witness}
\begin{tabular}{@{}rrrrrrrr@{\hspace{8mm}}rrrrrrrr@{}}\toprule
$i$&$j$&$k$&$D$&$Y$&$Z$&$T$&$m$&
$i$&$j$&$k$&$D$&$Y$&$Z$&$T$&$m$\\\midrule
1 & 2 & 11 & 7 & 1 & $-1$ & 4 & 5 & 3 & 7 & 11 & 24 & 2 & 0 & 9 & 7 \\
1 & 3 & 11 & 25 & 2 & 0 & 11 & 12 & 3 & 8 & 11 & 21 & 2 & 1 & 1 & 30 \\
1 & 4 & 10 & 8 & 1 & $-3$ & 1 & 4 & 3 & 9 & 7 & 3 & 0 & 0 & 1 & 1 \\
1 & 5 & 11 & 30 & $-1$ & 1 & 1 & 4 & 4 & 5 & 2 & 8 & 0 & 0 & 1 & 9 \\
1 & 6 & 9 & 7 & 0 & 2 & 3 & 3 & 4 & 6 & 7 & 50 & 1 & 1 & 18 & 4 \\
1 & 7 & 9 & 7 & 0 & 1 & 1 & 6 & 4 & 7 & 10 & 11 & 1 & 0 & 3 & 1 \\
1 & 8 & 10 & 17 & 1 & 1 & 2 & 4 & 4 & 8 & 2 & 6 & 1 & 0 & 1 & 20 \\
1 & 9 & 11 & 18 & 2 & $-1$ & 9 & 1 & 4 & 9 & 7 & 14 & 0 & 0 & 5 & 1 \\
2 & 3 & 7 & 12 & 0 & $-1$ & 1 & 2 & 5 & 6 & 7 & 7 & 1 & 0 & 1 & 23 \\
2 & 4 & 7 & 10 & $-1$ & 0 & 3 & 2 & 5 & 7 & 11 & 32 & $-1$ & 1 & 1 & 4 \\
2 & 5 & 7 & 9 & 0 & $-1$ & 1 & 8 & 5 & 8 & 11 & 8 & 1 & 0 & 1 & 22 \\
2 & 6 & 6 & 3 & 0 & 0 & 1 & 1 & 5 & 9 & 7 & 8 & 1 & $-1$ & 1 & 3 \\
2 & 7 & 11 & 16 & 2 & $-1$ & 8 & 8 & 6 & 7 & 1 & 26 & 3 & $-1$ & 12 & 17 \\
2 & 8 & 11 & 7 & 2 & $-1$ & 1 & 8 & 6 & 8 & 1 & 5 & 3 & $-1$ & 1 & 12 \\
2 & 9 & 6 & 3 & 0 & 0 & 2 & 1 & 6 & 9 & 1 & 3 & 0 & 0 & 1 & 1 \\
3 & 4 & 4 & 3 & 0 & 0 & 1 & 1 & 7 & 8 & 10 & 8 & 0 & 1 & 1 & 10 \\
3 & 5 & 3 & 10 & 1 & 0 & 1 & 5 & 7 & 9 & 11 & 17 & 2 & $-1$ & 8 & 10 \\
3 & 6 & 7 & 8 & 0 & 1 & 1 & 2 & 8 & 9 & 11 & 8 & 2 & $-1$ & 1 & 5 \\
\bottomrule\end{tabular}
\end{table}

For instance, $P_2$ and $P_6$ lie on opposite sides of $\{F_6=0\}$. At
$(115/3,0,0,1/3)$ we have $F_6=0$; the remaining signed forms of $P_2$
($-F_5,F_7,F_{11}$) and of $P_6$ ($F_1,F_7,F_9$) take the values
$21/3,\,4/3,\,59/3$ and $1/3,\,4/3,\,36/3$.

\begin{remark}[Contacts above the base]\label{rem:heights}
The nine bases $\beta_i$ are tetrahedra in the hyperplane $\{t=0\}\cong\R^3$.
They cannot be pairwise neighbourly, since $S_3=8$~\cite{Z91}. Indeed, an
exact computation shows that for $25$ pairs $(i,j)$ the set
$P_i\cap P_j\cap\{t=0\}=\beta_i\cap\beta_j$ is $2$-dimensional, while for the
$11$ pairs
\[
 (1,2),(1,3),(1,6),(1,9),(2,4),(2,7),(3,7),(4,6),(4,9),(6,7),(7,9)
\]
the bases are disjoint and $P_i$, $P_j$ touch only at positive heights. For
example, on $P_1\cap P_2$ we have $F_{11}=0$ and
$F_8,\,-F_9,\,-F_5,\,F_7\ge0$, and the identity
\[
   13F_8-14F_9-F_5+3F_7-F_{11}=195\,t-38
\]
shows that $t\ge38/195$ on $P_1\cap P_2$.
\end{remark}

The vertices of $P_1,\dots,P_9$ are listed in Appendix~\ref{app:vertices}.

\section{Doubling and the proof of \texorpdfstring{Theorem~\ref{thm:A}}{Theorem A}}\label{sec:doubling}
Write points of $\R^d$ as $(x_1,\dots,x_{d-1},t)$ and let $H=\{t=0\}$ and
$H^{\pm}=\{\pm t\ge0\}$.

\begin{definition}\label{def:based}
A family of $d$-simplices in $\R^d$ is \emph{based} if each member $P$ lies in
$H^+$ and has a facet $\beta(P)\subseteq H$, its \emph{base}. A based family
is \emph{pierced} if some line in $H$ meets the relative interior of every
base. An arbitrary family of $d$-simplices is \emph{pierced} if some line
meets the interior of every member.
\end{definition}

By Propositions~\ref{prop:simplices} and~\ref{prop:nine}, $P_1,\dots,P_9$
form a based, pierced, neighbourly family in $\R^4$; the piercing line is the
$x$-axis. We use two simple facts. First, if $P$ is a member of a based
family and $y\in\relint\beta(P)$, then $y+\tau e_d\in\interior P$ for all
sufficiently small $\tau>0$; indeed, all facet inequalities of $P$ other than
$t\ge0$ are strict at $y$. Second, affine bijections preserve
neighbourliness.

\begin{lemma}[Lifting]\label{lem:lift}
If there is a pierced neighbourly family of $m$ $d$-simplices in $\R^d$, then
there is a based, pierced, neighbourly family of $m$ $(d+1)$-simplices in
$\R^{d+1}$.
\end{lemma}

\begin{proof}
Identify $\R^d$ with $H\subset\R^{d+1}$, choose a point $a$ with positive last
coordinate, and replace each member $P$ by the pyramid $\conv(P\cup\{a\})$.
Every point of such a pyramid other than $a$ can be written uniquely as
$a+\lambda(p-a)$ with $p\in H$ and $0<\lambda\le1$, and then $p\in P$. Hence
$\conv(P\cup\{a\})\cap\conv(Q\cup\{a\})=\conv((P\cap Q)\cup\{a\})$, which is
$d$-dimensional when $P\cap Q$ is $(d-1)$-dimensional. The bases are the
original simplices, and their interiors in $H$ are met by the original line.
\end{proof}

\begin{lemma}[Doubling, Zaks~\cite{Z81}]\label{lem:double}
Let $d\ge3$. If there is a based, pierced, neighbourly family of $m$
$d$-simplices in $\R^d$, then there is a pierced neighbourly family of $2m$
$d$-simplices in $\R^d$.
\end{lemma}

\begin{figure}[ht]
\centering
\begin{tikzpicture}[scale=0.45]
  \definecolor{oiblue}{HTML}{0072B2}\definecolor{oired}{HTML}{D55E00}
  \draw[->,gray] (-4.2,0)--(11.6,0) node[right,black]{\small$x_1$};
  \draw[->,gray] (0,-7)--(0,4.8) node[above,black]{\small$x_2$};
  \foreach \a in {1,2,3}{
    \fill[oiblue,opacity=0.25] (\a-0.18,-6.5) rectangle (\a+0.18,4.4);
    \draw[oiblue] (\a-0.18,-6.5) rectangle (\a+0.18,4.4);
    \fill[oired,opacity=0.25] (-3.9,\a-0.18) rectangle (10.9,\a+0.18);
    \draw[oired] (-3.9,\a-0.18) rectangle (10.9,\a+0.18);}
  \foreach \a in {1,2,3}{\foreach \b in {1,2,3}{\fill (\a,\b) circle (2.2pt);}}
  \draw[thick,dashed] (0.2,-6.8)--(11.3,4.3);
  \foreach \a in {1,2,3}{\draw[thick] (\a,\a-7) circle (4pt);
     \draw[thick] (7+\a,\a) circle (4pt);}
  \node[oiblue,below] at (2,-6.6) {\small$\sigma(\beta_i)$};
  \node[oired,left] at (-3.9,2) {\small$\tau(\beta_j)$};
  \node[right] at (3.2,3.5) {\small$q_{ij}$};
  \draw[->,thin] (3.2,3.5) -- (3.08,3.1);
  \node[right] at (11.3,4.3) {\small$\ell$};
\end{tikzpicture}
\caption{The hyperplane $H$ in the proof of Lemma~\ref{lem:double}
(schematic, $d=3$). The stretched bases $\sigma(\beta_i)$ (blue; the
simplices $\sigma(P_i)$ lie above $H$) and $\tau(\beta_j)$ (orange;
$\tau(P_j)$ lies below $H$) overlap around the points $q_{ij}$. The
projection of the line $\ell$ (dashed) meets the blue bases at $s=a_i$ and
the orange bases at $s=c+a_j$ (circles); $\ell$ itself passes slightly above
$H$ for $s<c/2$ and slightly below for $s>c/2$.}\label{fig:doubling}
\end{figure}
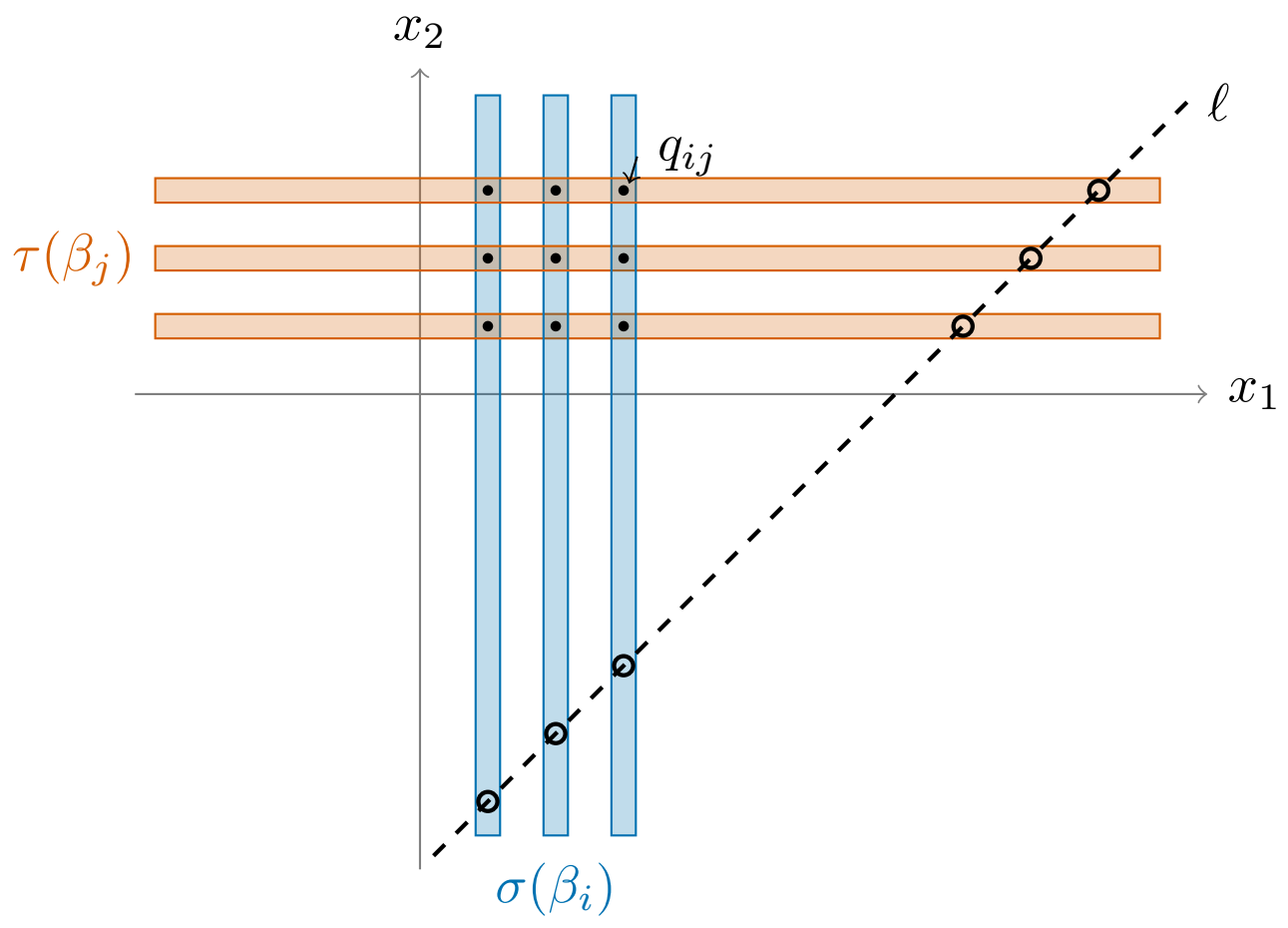

\begin{proof}
Let $P_1,\dots,P_m$ be the based family and $\beta_i=\beta(P_i)$. An affine map
$(x,t)\mapsto(Ax+b,t)$ preserves $H$ and $H^+$, so we may assume that the
piercing line is the $x_1$-axis. Choose $a_i$ with
$a_ie_1\in\relint\beta_i$ and $\varepsilon>0$ such that the ball of radius
$\varepsilon$ in $H$ centred at $a_ie_1$ lies in $\beta_i$ for every $i$. Let
$R=\max_i|a_i|$, fix $c>2R$, and choose $L>0$ with $L\varepsilon>c+R$.
Since $d-1\ge2$, the hyperplane $H$ contains the coordinates $x_1,x_2$; write
$x=(x_1,x_2,x')$ and define the affine bijections
\[
  \sigma(x_1,x_2,x',t)=(x_1,\,Lx_2,\,x',\,t),\qquad
  \tau(x_1,x_2,x',t)=(Lx_2,\,x_1,\,x',\,-t).
\]
We claim that $\{\sigma(P_i)\}_{i\le m}\cup\{\tau(P_j)\}_{j\le m}$ is the
required family. The simplices $\sigma(P_i)$ are pairwise neighbourly, and so
are the $\tau(P_j)$.

Since $\sigma(P_i)\subseteq H^+$ and $\tau(P_j)\subseteq H^-$ (see
Figure~\ref{fig:doubling}), we have
$\sigma(P_i)\cap\tau(P_j)=\sigma(\beta_i)\cap\tau(\beta_j)\subseteq H$. Let
$q=(a_i,a_j,0,\dots,0)\in H$. Then $\sigma^{-1}(q)=(a_i,a_j/L,0,\dots,0)$ and
$\tau^{-1}(q)=(a_j,a_i/L,0,\dots,0)$ are at distance less than $\varepsilon$
from $a_ie_1$ and $a_je_1$, respectively. Since $\sigma$ and $\tau$ restrict to
homeomorphisms of $H$, the point $q$ has a neighbourhood in $H$ contained in
$\sigma(\beta_i)\cap\tau(\beta_j)$. Thus $\sigma(P_i)$ and $\tau(P_j)$ are
neighbourly.

For $\delta>0$ consider the line
$\ell(s)=(s,\,s-c,\,0,\dots,0,\,\delta(c/2-s))$. For $s=a_i$ we get
\[
  \sigma^{-1}(\ell(a_i))=\bigl(a_i,\,(a_i-c)/L,\,0,\dots,0,\,
     \delta(c/2-a_i)\bigr),
\]
whose projection to $H$ is within $(c+R)/L<\varepsilon$ of $a_ie_1$, and whose
last coordinate is positive because $c/2>R$. For $s=c+a_j$ we get
\[
  \tau^{-1}(\ell(c+a_j))=\bigl(a_j,\,(c+a_j)/L,\,0,\dots,0,\,
     \delta(c/2+a_j)\bigr),
\]
again with projection within $\varepsilon$ of $a_je_1$ and with positive last
coordinate. By the first fact above, for $\delta$ small enough these $2m$
points lie in the interiors of the corresponding simplices $P_i$ and $P_j$. Hence $\ell$ meets the
interiors of all $2m$ simplices.
\end{proof}

\begin{proof}[Proof of Theorem~\ref{thm:A}]
Lemma~\ref{lem:double} applied to $P_1,\dots,P_9$ gives $18$ pairwise
neighbourly $4$-simplices in $\R^4$, together with a line meeting all their
interiors. If $\cF$ is a pierced neighbourly family of $m$ simplices in
$\R^d$, $d\ge4$, then Lemmas~\ref{lem:lift} and~\ref{lem:double} give a
pierced neighbourly family of $2m$ simplices in $\R^{d+1}$. By induction,
$S_d\ge18\cdot2^{d-4}$ for every $d\ge4$.
\end{proof}

\begin{corollary}\label{cor:based3}
A based, pierced, neighbourly family of tetrahedra in $\R^3$ has at most four
members.
\end{corollary}

\begin{proof}
By Lemma~\ref{lem:double}, five such tetrahedra would give ten pairwise
neighbourly tetrahedra, contrary to $S_3=8$~\cite{Z91}.
\end{proof}

\begin{remark}[An explicit family]\label{rem:explicit}
The proof of Lemma~\ref{lem:double} can be carried out explicitly for
$P_1,\dots,P_9$. Let $L=100\,000$ and
\[
  A(x,y,z,t)=(x,\,Ly,\,z,\,t),\qquad B(x,y,z,t)=(Ly,\,x,\,z,\,-t).
\]
The $18$ simplices $A_i=A(P_i)$ and $B_j=B(P_j)$ are pairwise neighbourly.
The point $Q_{ij}=(u_i,u_j,0,0)$ lies in the relative interiors of the bases
of $A_i$ and of $B_j$: its preimages $(u_i,u_j/L,0,0)\in\beta_i$ and
$(u_j,u_i/L,0,0)\in\beta_j$ differ from $(u_i,0,0,0)$ and $(u_j,0,0,0)$ only in
the $y$-coordinate, by at most $40/L$. By Proposition~\ref{prop:simplices}
each signed side form is at least $1/2$ at $(u_i,0,0,0)$, and it changes by at
most $212\cdot40/L<1/2$. Similarly the line
\[
   \ell(s)=\bigl(s,\,s-64,\,0,\,(52-s)/L\bigr)
\]
meets the interior of $A_i$ at $s=u_i$ and of $B_j$ at $s=64+u_j$. The
preimage of the point in $P_i$ (respectively $P_j$) has first coordinate $u_i$
(respectively $u_j$), third coordinate $0$, positive last coordinate at most
$52/L$, and second coordinate of absolute value at most $104/L$. The signed
side forms therefore change by at most $(212\cdot104+103\cdot52)/L<1/2$. The
vertices of all $18$ simplices are obtained from Appendix~\ref{app:vertices}
by applying $A$ and $B$.

Applying Lemmas~\ref{lem:lift} and~\ref{lem:double} once more, with explicit
parameters, we have also verified in exact arithmetic a family of $36$
pairwise neighbourly $5$-simplices in $\R^5$ (all $630$ pairwise
intersections are $4$-dimensional); this exceeds $2^5=32$.
\end{remark}

\begin{remark}[Based families and codes]\label{rem:basedcodes}
In a based family all members share the facet hyperplane $H$, on the same
side. Deleting this common position from the code of Lemma~\ref{lem:coding}
leaves a neighbourly $d$-code. Thus a based neighbourly family in $\R^d$ has
at most $M_d$ members. For the family $P_1,\dots,P_9$ this code is exactly
$W$ (Table~\ref{tab:words}), and the code of the $18$ simplices of
Remark~\ref{rem:explicit} is $D(W)$. Since $M_4=12$
(Theorem~\ref{thm:B}), doubling a based family in $\R^4$ can never give more
than $24$ simplices. This matches the upper bound $S_4\le24$ of
Theorem~\ref{thm:C}.
\end{remark}

\section{Rainbow triples}\label{sec:rainbow}
Throughout this section $V$ is a neighbourly $k$-code. The \emph{rainbow
hypergraph} of $V$ has vertex set $V$, and its edges are the rainbow triples.
It is \emph{linear} if every pair of words lies in at most one rainbow
triple.

\subsection{Codes without rainbow triples}
\begin{lemma}[Splitting]\label{lem:split}
Let $A$ be a neighbourly $k$-code with at least two words and no rainbow
triple. Then some column contains both bits and has no star in $A$.
Consequently every neighbourly $k$-code without rainbow triples has at most
$2^{k-1}$ words, for $k\ge1$.
\end{lemma}

\begin{proof}
Among the columns containing both bits (they exist because any two words
conflict), choose a column $i$ with the largest number of bits. Let $A_0,A_1$
be its two sides. Suppose some $w\in A$ has $w_i=\st$. For $a\in A_0$ and
$b\in A_1$ we have $c(a,b)=i$, while $c(w,a)\ne i\ne c(w,b)$. Since
$\{w,a,b\}$ is not rainbow, $c(w,a)=c(w,b)$. Fixing $a$ and varying $b$, and
then fixing $b$ and varying $a$, we see that $c(w,\cdot)$ is constant on
$A_0\cup A_1$, say equal to $j$. Column $j$ then contains both bits and has
bits in all words of $A_0\cup A_1\cup\{w\}$, contradicting the choice of $i$.
Hence $A=A_0\cup A_1$.

For the bound we use induction on $k$. Two words of weight one with different
supports cannot conflict, so the bound holds for $k=1$. For $k\ge2$ and
$|A|\ge2$, delete column $i$ from both sides. Each side becomes a neighbourly
$(k-1)$-code without rainbow triples, because all conflicts inside a side lie
outside column $i$. Hence $|A_0|,|A_1|\le2^{k-2}$.
\end{proof}

\subsection{The neighbourhood of a rainbow triple}
\begin{lemma}\label{lem:neighbourhood}
Let the rainbow hypergraph of $V$ be linear, let $T_0=\{a,b,c\}$ be a rainbow
triple, and let $I$ be the set of its three colours. Then every word
$w\notin T_0$ is of one of two kinds.
\begin{enumerate}[(i)]
\item \emph{First kind:} $w$ conflicts with all three words of $T_0$ at the
same position $p_w\notin I$, and $w$ has stars on $I$.
\item \emph{Second kind:} $w$ conflicts with the three words of $T_0$ at
three different positions, namely the positions of $I$, and all words of the
second kind have the same bits on $I$.
\end{enumerate}
Let $X$ and $Y$ be the sets of words of the second and the first kind,
respectively. Then every $y\in Y$ conflicts at $p_y$ with every word of
$T_0\cup X$.
\end{lemma}

\begin{proof}
Suppose that $c(w,a)=c(w,b)=p$ and $c(w,c)=q\ne p$. Then $a,b$ have the same
bit at $p$. If $c(a,c)=p$, then also $c(b,c)=p$, which is impossible in a
rainbow triple; so neither $ac$ nor $bc$ has colour $p$. The triples
$\{w,a,c\}$ and $\{w,b,c\}$ share the pairs $ac$ and $bc$ with $T_0$, so by
linearity they are not rainbow. This forces $c(a,c)=q=c(b,c)$, a
contradiction. Hence the three colours $c(w,a),c(w,b),c(w,c)$ are either
equal or pairwise distinct.

If they are all equal to $p_w$, then the words of $T_0$ have the same bit at
$p_w$, so $p_w\notin I$; and since every position of $I$ carries both bits in
$T_0$, a bit of $w$ on $I$ would create a second conflict.

Suppose that they are pairwise distinct. For each pair $u,v$ of $T_0$, the
colour $c(u,v)$ equals one of $c(w,u),c(w,v)$; otherwise $\{w,u,v\}$ would be
a second rainbow triple containing $uv$. So $w$ has bits on all three
positions of $I$. At the position $c(u,v)$ only $u$ and $v$ have bits, so the
bit of $w$ there conflicts with exactly one of them. Each word of $T_0$ must be
hit exactly once. Orienting every edge of the triangle $T_0$ towards the word
it hits, we get one of the two cyclic orientations. The two corresponding
patterns on $I$ are complementary, so they cannot both occur (they would
conflict three times).

Finally let $x\in X$, $y\in Y$, and suppose that $c(x,y)\ne p_y$. The colours
$c(x,u)$, $u\in T_0$, are the three elements of $I$, so for at least two
$u\in T_0$ the colour $c(x,u)$ differs from $c(x,y)$; it also differs from
$c(y,u)=p_y\notin I$. Then $\{x,y,u\}$ is rainbow for two different $u$,
contrary to linearity.
\end{proof}

\subsection{Codes with a linear rainbow hypergraph}
A \emph{full partition} of weight $j$ is a set of $2^j$ words of weight $j$,
every two of which conflict exactly once; here equal supports are allowed.
The cylinders of a full partition tile the cube.

\begin{lemma}\label{lem:fullpartition}
Every full partition of weight $j$ is, up to isomorphism, the set of
root-to-leaf paths of a complete binary tree of depth $j$ whose internal
nodes carry pairwise distinct positions: the word of a leaf records the bits
of the branches taken along its path.
\end{lemma}

\begin{proof}
Let $j>0$ and let $S$ be the support of some word. The character
$\chi_S(x)=(-1)^{\sum_{i\in S}x_i}$ has mean zero on the cube and on every
cylinder whose support differs from $S$ (such a support misses some
position of $S$), and it is constant on every cylinder with support $S$.
Since the cylinders tile the cube and have equal measure, the words with
support $S$ come in pairs with opposite values of $\chi_S$. Three words with
the same support cannot pairwise conflict exactly once, so every word has
exactly one \emph{twin}, with the same support and a single conflict. At the
conflict position of two twins every other word has a star, since otherwise
its numbers of conflicts with the two twins would differ by one. Merging
every pair of twins into one word with a star at that position gives a full
partition of weight $j-1$, and induction finishes the proof.
\end{proof}

\begin{theorem}\label{thm:linear}
Let $k\ge2$, and let $V$ be a neighbourly $k$-code with a linear rainbow
hypergraph. Then the rainbow triples of $V$ are pairwise disjoint and
$|V|\le3\cdot2^{k-2}$. Equality holds if and only if $V$ is isomorphic to
$G_k$.
\end{theorem}

\begin{proof}
\emph{Disjointness.} By linearity two rainbow triples share at most one
word. Let $T_0$ be a rainbow triple and let $u,v\notin T_0$; we show that
$\{u,v,x\}$ is not rainbow for $x\in T_0$, using
Lemma~\ref{lem:neighbourhood}. If $u,v\in Y$, then $c(u,x)=p_u$ and
$c(v,x)=p_v$ do not depend on $x$, so if $\{u,v,x\}$ were rainbow for one
$x$, it would be rainbow for all three, which contradicts linearity. If
$u,v\in X$, then $c(u,x)=c(v,x)$ for every $x\in T_0$. If $u\in Y$ and $v\in
X$, then $c(u,v)=p_u=c(u,x)$.

\emph{Bound.} Let $m$ be the number of rainbow triples and $n=|V|$. Deleting
one word from each triple leaves a code without rainbow triples, so
$n-m\le2^{k-1}$ by Lemma~\ref{lem:split}. Since $m\le n/3$, this gives
$2n/3\le2^{k-1}$.

\emph{Equality.} Let $n=3\cdot2^{k-2}$. Then $m=2^{k-2}$ and the rainbow
triples cover $V$. Choosing two words from each triple, we get a code $A$
without rainbow triples and with $2^{k-1}$ words. In such a code every mixed
column is \emph{balanced}, i.e.\ has as many zeros as ones. Indeed, the
column $i$ of Lemma~\ref{lem:split} splits $A$ into two sides of size
exactly $2^{k-2}$. Any other mixed column has bits in one side only: words
from different sides already conflict at $i$, so their bits elsewhere must
agree. We conclude by induction applied to that side.

Suppose that some word $w$ is of the second kind with respect to a rainbow
triple $T$. Pick $i$ in the colour set of $T$. Form $A$ by keeping the two
words of $T$ with bits at $i$, the word $w$ with any partner from its own
triple, and two arbitrary words from every other triple. By
Lemma~\ref{lem:neighbourhood}, every word of $A\setminus T$ with a bit at $i$
has the same bit as $w$, so column $i$ of $A$ is not balanced, a
contradiction. Hence, with respect to every rainbow triple $T$, all words
outside $T$ are of the first kind. For two triples $T,T'$ this means that all
nine pairs between them have a common colour $p_{T,T'}$, at which all words
of $T$ have one bit and all words of $T'$ the other.

Let $C_T$ be the pattern of bits common to the three words of $T$. These
$2^{k-2}$ patterns pairwise conflict, at the positions $p_{T,T'}$, and each
has weight at most $k-2$, because every word of $T$ has two bits among the
colours of $T$. Their cylinders are disjoint and each has measure at least
$2^{-(k-2)}$, so all weights are exactly $k-2$ and the patterns form a full
partition. By Lemma~\ref{lem:fullpartition} it is the full binary tree of
depth $k-2$. Each triple $T$ consists of its core $C_T$ followed by a copy of
$T$ on its three colours, where all other words have stars. This is $G_k$.
Conversely, $G_k$ has a linear rainbow hypergraph: a triple meeting both
sides of a doubling has two conflicts at the root position, so the only
rainbow triples are the $2^{k-2}$ final copies of $T$.
\end{proof}

\subsection{Excluding triangles}
\begin{theorem}\label{thm:linearnocore}
Let $k\ge2$, and let $V$ be a neighbourly $k$-code with a linear rainbow
hypergraph and without triangles (Definition~\ref{def:triangle}). Then
$|V|\le2^{k-1}$.
\end{theorem}

\begin{proof}
We describe a replacement that preserves $|V|$, keeps the hypotheses, and
strictly decreases the number of rainbow triples. After finitely many
replacements no rainbow triple is left, and Lemma~\ref{lem:split} gives the
bound.

Let $T_0$ be a rainbow triple, with $X,Y$ as in
Lemma~\ref{lem:neighbourhood}, and put $B=T_0\cup X$; then $V=B\cup Y$. Let
$C$ be the pattern of bits common to all words of $B$, let $h$ be its weight
and $s=k-h$. No position of $I$ belongs to $C$, and every word of $T_0$ has
two bits on $I$, so $h\le k-2$. Equality would make $T_0$ a triangle, so
$s\ge3$.

\emph{How many words must be replaced?} All words of $X$ have the bits of $C$
and the same three bits on $I$. After these $h+3$ positions are deleted, they
have $s-3$ bits each, and every two of them still conflict. Their cylinders
are disjoint, so $|X|\le2^{s-3}$ and
\begin{equation}\label{eq:Bsize}
   |B|\le3+2^{s-3}\le2^{s-1}.
\end{equation}

\emph{A rainbow-free replacement.} On fresh positions take a complete binary
tree of depth $s-1$ with distinct positions at its internal nodes. To each
leaf assign the word that records the $s-1$ branch bits along its path,
together with one bit at a fresh private position of that leaf. Two such
words conflict only at the node where their paths split, so among any three
of them one splits off from the other two at a single node, and two of the
three colours coincide. By~\eqref{eq:Bsize} we can choose $|B|$ leaves.
Prefix each chosen word with $C$, and replace $B$ by the resulting words.
They have weight $h+s=k$ and distinct supports.

By Lemma~\ref{lem:neighbourhood}, every $y\in Y$ conflicts with every word of
$B$ at $p_y$, and only there. Thus $p_y$ is a position of $C$, and $y$
conflicts with $C$ exactly once. Since $y$ has stars on the fresh positions,
it keeps colour $p_y$ with every new word. So the new code is a neighbourly
$k$-code of the same size.

\emph{Rainbow triples of the new code.} Three new words do not form a rainbow
triple. Two new words and an outer word $y$ do not either, since both pairs
with $y$ have colour $p_y$. If a new word and two outer words $y,y'$ formed a
rainbow triple, then the colours $p_y$, $p_{y'}$, $c(y,y')$ would be
distinct. Then $\{b,y,y'\}$ would have been rainbow for every $b\in B$, so
the pair $yy'$ would have been in $|B|\ge3$ rainbow triples. Hence every
rainbow triple of the new code consists of old words of $Y$ and was a rainbow
triple before. The new hypergraph is therefore linear, it has lost at least
$T_0$, and it contains no triangle, since every triangle is a rainbow triple.
\end{proof}

\subsection{Back to geometry}
\begin{lemma}\label{lem:line}
The code of a neighbourly family of $d$-simplices in $\R^d$, $d\ge2$, contains
no triangle, i.e.\ no three words sharing $d-1$ bits.
\end{lemma}

\begin{proof}
Suppose that three members have codes forming a triangle, with core of
weight $d-1$ and local positions $p,q,r$ as in Section~\ref{sec:codes}. The
$d-1$ core hyperplanes span facets of each of the three simplices. Any $d-1$
facets of a $d$-simplex meet in an edge, so these hyperplanes meet in a line
$\Lambda$, and each of the three simplices meets $\Lambda$ in one of its
edges. The ends of these edges lie on the hyperplanes of the remaining two
facets. Let $P,Q,R$ be the points where $\Lambda$ meets the hyperplanes with
indices $p,q,r$. Up to isomorphism the three words read
$00\st,\,1\st0,\,\st11$ on $(p,q,r)$, so the three edges are $[P,Q]$,
$[P,R]$ and $[Q,R]$, and the points $P,Q,R$ are pairwise distinct. The first
two simplices conflict at $p$, so they lie on opposite sides of the
hyperplane through $P$. Hence $Q$ and $R$ lie on opposite sides of $P$ on
$\Lambda$, i.e.\ $P$ lies strictly between $Q$ and $R$. The first and third
simplices conflict at $q$, and the same argument places $Q$ strictly between
$P$ and $R$. Both cannot hold on a line.
\end{proof}

\begin{proof}[Proof of Theorem~\ref{thm:D}]
By Lemma~\ref{lem:coding}, rainbow triples of $\cF$ correspond to rainbow
triples of its code $w(\cF)$, a neighbourly $(d+1)$-code, and by
Lemma~\ref{lem:line} the code contains no triangle. If the rainbow
hypergraph is linear, Theorem~\ref{thm:linearnocore} with $k=d+1$ gives
$|\cF|\le2^d$. In general, take a maximal collection of pairwise disjoint
quadruples $\{A,B,C,D\}$ with $\{A,B,C\}$ and $\{A,B,D\}$ rainbow, say $q$ of
them. Removing their $4q$ members leaves a family with a linear rainbow
hypergraph. Removing members does not change the colours of the remaining
pairs, and a remaining pair in two rainbow triples would give another
quadruple. Hence $|\cF|-4q\le2^d$, and so $q\ge u/4$.
\end{proof}

\begin{example}\label{ex:rainbow18}
(a) Such quadruples occur naturally. Let $O$ be an interior point of a
tetrahedron and join $O$ to its four facets. The resulting four tetrahedra
are pairwise neighbourly, and the contact of the two cones over facets $F_i$
and $F_j$ is the triangle spanned by $O$ and the edge $F_i\cap F_j$. The three
edges involved in any three of the cones share a vertex, any two of them span
a facet plane, and $O$ is not on it. Hence every triple of the four tetrahedra
is rainbow.

(b) For the $18$ simplices of Remark~\ref{rem:explicit}, the code
$D(W)$ has $23$ positions and no triangle (as it must, by
Lemma~\ref{lem:line}). Every rainbow triple lies in one layer, since a triple
meeting both layers has two conflicts at the base position. There are $26$
rainbow triples in each layer, and $23$ pairs of each layer lie in at least
two of them. For instance, the pair $\{A_1,A_4\}$ lies in the five rainbow
triples $\{A_1,A_4,A_j\}$, $j=2,3,5,6,9$.
\end{example}

The upper bound $S_4\le24$ of Theorem~\ref{thm:C} would improve by one if
the extremal $5$-code were unique.

\begin{proposition}\label{prop:S4}
If every neighbourly $5$-code with $24$ words is isomorphic to $G_5$, then
$S_4\le23$.
\end{proposition}

\begin{proof}
A neighbourly family of $24$ simplices in $\R^4$ would have a neighbourly
$5$-code of size $24=M_5$, which would then be isomorphic to $G_5$. But $G_5$
contains triangles (the final copies of $T$ with the three common bits added
by the doublings), and Lemma~\ref{lem:line} excludes them.
\end{proof}

By Theorem~\ref{thm:linear}, the hypothesis of Proposition~\ref{prop:S4} holds
for all $5$-codes of size $24$ with a linear rainbow hypergraph.

\section{Upper bounds for neighbourly codes}\label{sec:upper}
In this section $k\ge2$, volumes are measured in units of $2^{-k}$, and a
triangle in a neighbourly $k$-code always has a core of weight $k-2$.

\subsection{Isolated triangles}
\begin{lemma}[Isolation]\label{lem:isolation}
Let $V$ be a neighbourly $k$-code with $|V|=M_k$, and let $\mathcal T\subseteq V$
be a triangle with core $C$ and local positions $p,q,r$. Then every
$w\in V\setminus\mathcal T$ conflicts with $C$ exactly once and has stars at
$p,q,r$.
\end{lemma}

\begin{proof}
Normalise $\mathcal T$ to $00\st,\ 1\st0,\ \st11$ on $(p,q,r)$. If $w$
conflicts with $C$, this conflict is its conflict with each word of
$\mathcal T$. A bit of $w$ at $p$, $q$ or $r$ would create a second conflict
with one of the words, since each of these columns contains both bits in
$\mathcal T$.

Suppose that $w$ does not conflict with $C$. A single bit at one of $p,q,r$
conflicts with exactly one word of $\mathcal T$, so $w$ must have bits at all
three positions; this is impossible for $k=2$. For $k\ge3$ let $x$ be the bit
of $w$ at $p$. A single conflict with $00\st$ forces $w_q=1-x$, and a single
conflict with $1\st0$ forces $w_r=x$. So $w$ reads $010$ or $101$ on
$(p,q,r)$. These two patterns conflict three times, so every word outside
$\mathcal T$ has on $\{p,q,r\}$ either three stars or the same pattern as
$w$. Hence on these positions $w$ conflicts only with words of $\mathcal T$.

Now delete $w$ and $\mathcal T$. Let $P$ be the pattern of the $k-3$ bits of $w$
outside $\{p,q,r\}$, and add the six words obtained by prefixing $P$ to the
words of $D(T)$, placed on seven new positions. Every remaining word conflicts
exactly once with $P$, and the new words have distinct supports and conflict
pairwise once. The result is a neighbourly $k$-code with $M_k-4+6$ words, a
contradiction.
\end{proof}

\begin{corollary}\label{cor:cores}
Let $|V|=M_k$. Then distinct triangles of $V$ are disjoint, and the cylinders
of the cores of all triangles together with the cylinders of the words
outside triangles are pairwise disjoint. In particular $|V|+t\le2^k$, where
$t$ is the number of triangles.
\end{corollary}

\begin{proof}
Two words of a triangle share $k-1$ positions of their supports, while by
Lemma~\ref{lem:isolation} a word outside a triangle shares at most $k-2$ with
each of its words. Hence triangles are disjoint. A word outside a triangle
conflicts with its core. For two triangles, the conflict between any of their
words lies, by Lemma~\ref{lem:isolation}, in both cores. The core cylinders
have volume $4$, so $(|V|-3t)+4t\le2^k$.
\end{proof}

\subsection{Exchange}\label{subsec:exchange}
A neighbourly $k$-code $V$ is \emph{normalised} if $|V|=M_k$ and $V$ has the
largest number $t$ of triangles among codes of size $M_k$. For $v\in V$ and
$i\in\supp v$ let
\[
   c_i(v)=|V^{i,1-v_i}|
\]
be the size of the side opposite to $v$ in column $i$. Each word other than
$v$ conflicts with $v$ exactly once, so $\sum_{i\in\supp v}c_i(v)=|V|-1$. Let
$c_1\le\dots\le c_k$ denote the values $c_i(v)$ in nondecreasing order.

\begin{lemma}[Exchange]\label{lem:exchange}
Let $V$ be normalised, $v\in V$ and $J\subseteq\supp v$ with $|J|=r\ge2$.
Then $\sum_{i\in J}c_i(v)\ge3\cdot2^{r-2}-1$. If $v$ lies in no triangle,
then $\sum_{i\in J}c_i(v)\ge3\cdot2^{r-2}$. Consequently
\begin{equation}\label{eq:prefix}
  c_1+\dots+c_j\ \ge\ 3\cdot2^{j-2}-1\quad(2\le j\le k),
\end{equation}
and the bound improves to $3\cdot2^{j-2}$ if $v$ lies in no triangle.
\end{lemma}

\begin{proof}
Delete $v$ and all words conflicting with $v$ at a position of $J$; these are
$1+\sum_{i\in J}c_i(v)$ words. Let $P$ be the pattern of the $k-r$ bits of $v$
outside $J$, and add the $3\cdot2^{r-2}$ words of $G_r$, placed on new
positions and prefixed by $P$. Every remaining word conflicts with $v$ outside
$J$, hence exactly once with $P$; so we get a neighbourly $k$-code. If fewer
than $3\cdot2^{r-2}$ words were deleted, its size exceeds $M_k$. This proves
the first bound.

Let $v$ lie in no triangle, and suppose that exactly $3\cdot2^{r-2}$ words
were deleted. By Lemma~\ref{lem:isolation}, $v$ conflicts with all three
words of an old triangle at the same position, which lies in that triangle's
core. So every old triangle was deleted entirely or not at all. Old triangles
are disjoint and do not contain $v$, so at most
$(3\cdot2^{r-2}-1)/3<2^{r-2}$ of them were deleted. The $2^{r-2}$ triangles
of $G_r$, prefixed by $P$, are triangles of the new code. Thus the new code
has size $M_k$ and more triangles than $V$, a contradiction.
\end{proof}

\subsection{Private pairs and a rank inequality}\label{subsec:rank}
From now on we assume:
\begin{center}
\emph{$k\ge3$, $g=g_k$, $n=M_k>g$, and $V$ is a normalised code with $t$
triangles; $\Delta=2^k-n$.}
\end{center}

\noindent
By~\eqref{eq:prefix}, a word of a triangle has $c_i=1$ at its two local
positions (whose columns have profile $(1,1)$ by Lemma~\ref{lem:isolation})
and $c_i\ge3$ at the $k-2$ core positions. A word outside triangles has
$c_1+c_2\ge3$, so at most one of its columns has an opposite side of size at
most one.

Two words outside triangles that conflict in a column of profile $(1,1)$ form
a \emph{private pair}, and that column is a \emph{private column}. A column
of profile $(1,1)$ that contains a word of a triangle is a local column of
that triangle. Indeed, core columns have at least three words on one side.
So every $(1,1)$ column is either local or private, and private pairs are
disjoint. Let $m$ be the number of private pairs.

For words outside triangles let $n_0$ and $n_2$ be the numbers of those with
$c_1=0$ and with $c_1\ge2$, respectively. A word outside triangles with
$c_1=1$ either belongs to a private pair or lies on the larger side of a
column of profile $(a,1)$ with $a\ge2$; let $s$ be the number of words of the
second type. Then
\begin{equation}\label{eq:partition}
   n=3t+2m+s+n_2+n_0 .
\end{equation}
Let $h$ be the number of mixed columns of $V$, and for $v\in V$ put
\[
  F(v)=\frac12\sum_{i\in\supp v}\frac1{\max\{1,c_i(v)\}} .
\]

\begin{lemma}\label{lem:reciprocal}
Let $0<z_1\le\dots\le z_k$ and $0<y_1\le\dots\le y_k$ satisfy
$\sum_{i\le j}z_i\ge\sum_{i\le j}y_i$ for every $j$. Then
$\sum_i1/z_i\le\sum_i1/y_i$.
\end{lemma}

\begin{proof}
For a nondecreasing sequence and $u>0$ we have
$\sum_i(u-z_i)_+=\max_{0\le j\le k}\bigl(ju-\sum_{i\le j}z_i\bigr)$, so the
left-hand side for $z$ is at most that for $y$. Multiply by $2u^{-3}$ and
integrate over $u>0$, using $\int_z^\infty2(u-z)u^{-3}\,du=1/z$.
\end{proof}

\begin{table}[ht]
\centering\small
\caption{Comparison sequences and upper bounds for $F(v)$. Each sequence
continues with $6,12,24,\dots,3\cdot2^{k-3}$.}\label{tab:comparison}
\begin{tabular}{@{}lcc@{}}
\toprule
Class of $v$ & first three terms & upper bound for $F(v)$\\
\midrule
outside triangles, $c_1=0$ & $1,3,3$ & $1-1/g$\\
outside triangles, $c_1=1$ & $1,2,3$ & $13/12-1/g$\\
outside triangles, $c_1\ge2$ & $2,2,2$ & $11/12-1/g$\\
in a triangle & $1,1,3$ & $4/3-1/g$\\
\bottomrule
\end{tabular}
\end{table}

\begin{lemma}\label{lem:Fsum}
We have
\begin{equation}\label{eq:cost}
  h+\frac{n_0}2+\frac s8\ \le\ \sum_{v\in V}F(v)\ \le\
  n\Bigl(1-\frac1g\Bigr)+t+\frac{2m+s-n_2}{12}.
\end{equation}
\end{lemma}

\begin{proof}
For each class in Table~\ref{tab:comparison}, the prefix sums of
$(\max\{1,c_i\})_i$ dominate those of the comparison sequence $y$, by
Lemma~\ref{lem:exchange}. For $c_1=0$, replacing $0$ by $1$ adds one to every
prefix sum; for words in a triangle we use the values $1,1$ and $c_3\ge3$
noted above. Lemma~\ref{lem:reciprocal} gives
$F(v)\le\frac12\sum1/y_i$. Since
$\sum_{i\ge4}1/y_i=\frac13-\frac2g$, this is the bound in the table. Summing
over $V$ with~\eqref{eq:partition} gives the upper bound.

For the lower bound we sum by columns. A mixed column with sides of sizes
$a,b$ contributes $\frac12(a/b+b/a)\ge1$. If $b=1$ and $a\ge2$, the excess is
$(a-1)^2/(2a)\ge a/8$, and these columns have $s$ words on their larger sides
in total. A column with one empty side contributes $\frac12$ for each of its
bits, and such bits belong exactly to the $n_0$ words with $c_1=0$ (each of
which has $c_2\ge3$).
\end{proof}

Next we \emph{contract} each triangle to its core. This removes the $3t$
local columns and keeps all other columns mixed, and the resulting
$n-2t$ words still pairwise conflict exactly once. Remove also the $m$ private
columns. Two words are now joined if they conflict in a remaining column.
The resulting graph is complete multipartite with $N=n-2t-m$ parts: the
private pairs and singletons. Let $h'=h-3t-m$ be the number of remaining
mixed columns. Let $B$ be the $(n-2t)\times2h'$ matrix whose two columns for
each remaining mixed column are the indicator vectors of its $0$-side and
its $1$-side. For each private pair $\{u,v\}$ form the difference $B_u-B_v$
of rows. Let $r$ be the dimension of the space spanned by these $m$
differences.

\begin{lemma}[Rank inequality]\label{lem:rank}
$h\ge n+t-1+r$. Moreover, every difference $B_u-B_v$ is nonzero; in
particular $r\ge1$ if $m\ge1$.
\end{lemma}

\begin{proof}
Let $J$ be the direct sum of $h'$ blocks
$\bigl(\begin{smallmatrix}0&1\\1&0\end{smallmatrix}\bigr)$. Then $BJB^{\mathsf
T}$ is the adjacency matrix of the complete multipartite graph, and its
quadratic form is $(\sum_iX_i)^2-\sum_iX_i^2$, where $X_i$ is the sum of the
entries of the vector over part $i$. Choose one representative in each part,
and let $U$ be the space of vectors supported on the representatives with
$\sum_iX_i=0$. Then $\dim U=N-1$ and the form is negative definite on $U$.
Hence $B^{\mathsf T}$ is injective on $U$, and the form $J$ is negative
definite on $W_-=B^{\mathsf T}U$.

The two ends of a private pair have equal rows in the adjacency matrix. So
$z^{\mathsf T}BJ(B_u-B_v)^{\mathsf T}=0$ for every $z$, i.e.\ each difference
is $J$-orthogonal to the whole image of $B^{\mathsf T}$, which contains all the
differences. Hence the differences span a totally isotropic subspace of the
$J$-orthogonal complement of $W_-$. On this complement the form has negative
index $h'-(N-1)$, and a totally isotropic subspace has dimension at most the
negative index. Therefore $r\le h'-N+1$, which is the claim.

For nonvanishing, recall that the ends $u,v$ of a private pair lie outside
triangles. So all their bits outside the private column lie in mixed,
non-private, non-local columns (because $c_1+c_2\ge3$). These columns are
represented in $B$. Since $\supp u\ne\supp v$, some such position is in one
support only, and there $B_u-B_v$ has a nonzero entry.
\end{proof}

Delete the private bit from both words of every private pair; the two
resulting words of weight $k-1$ no longer conflict. Let $\ell_j\ge1$ be the
number of positions in the support of the first of them but not of the
second (and vice versa).

\begin{lemma}[Volume]\label{lem:volpairs}
The enlarged cylinders of the private pairs, the cores of the triangles and
the cylinders of the remaining words are pairwise disjoint. Consequently
\begin{equation}\label{eq:volpairs}
   \Delta\ \ge\ t+\sum_{j=1}^m\bigl(2-2^{1-\ell_j}\bigr)\ \ge\ t+m .
\end{equation}
\end{lemma}

\begin{proof}
Two words from different groups conflict outside the private columns, which
contain only their own pair; together with Corollary~\ref{cor:cores} this
gives disjointness. The two enlarged cylinders of the $j$-th pair have volume
$2$ each and intersect in a cylinder fixing $k-1+\ell_j$ bits, of volume
$2^{1-\ell_j}$. So their union has volume $4-2^{1-\ell_j}$.
\end{proof}

\begin{proposition}\label{prop:firstbalance}
Under the standing assumptions $r\ge1$ and
\begin{align}
  2m&\ \ge\ n_2+6n_0+\tfrac s2+12\bigl(\tfrac ng-1+r\bigr),\label{eq:balance1}\\
  n+24\bigl(\tfrac ng-1+r\bigr)&\ \le\ 6m+3t\ \le\ 6\Delta .\label{eq:balance2}
\end{align}
In particular $n(7+24/g)\le6\cdot2^k$.
\end{proposition}

\begin{proof}
Inequality~\eqref{eq:balance1} follows from~\eqref{eq:cost} and
Lemma~\ref{lem:rank}. Its right-hand side is positive because $n>g$, so
$m\ge1$ and $r\ge1$. Doubling~\eqref{eq:balance1}, adding $2m+3t$ and
using~\eqref{eq:partition} gives
$6m+3t\ge n+n_2+11n_0+24(n/g-1+r)$. The upper bound
in~\eqref{eq:balance2} is~\eqref{eq:volpairs}. Finally use $r\ge1$ and
$\Delta=2^k-n$.
\end{proof}

\begin{corollary}\label{cor:M5}
We have $M_5=24$ and $S_4\le24$. Moreover $M_k\le6\cdot2^k/(7+24/g_k)$
for every $k\ge5$; in particular $\limsup_kM_k/2^k\le6/7$.
\end{corollary}

\begin{proof}
For $k=5$ we have $g=24$, and $n>24$ would give $8n\le192$. For $k\ge5$ we
have $g_k(7+24/g_k)=7g_k+24\le8g_k=6\cdot2^k$, so the bound also covers the
case $M_k=g_k$. The last claim follows from Lemma~\ref{lem:coding}.
\end{proof}

\subsection{Columns with small sides}\label{subsec:small}
We keep the standing assumptions and assume moreover $k\ge4$. Let $Q=2^k$.
Since $n>g=3Q/4$ we have $n>3\Delta$, and it suffices to prove
$(46+78/g)n\le35Q$.

Let $e$ be the number of columns of profile $(2,1)$, and $H$ the total number
of words on the larger sides of columns of profile $(a,1)$ with $a\ge3$; thus
$s=2e+H$. These $2e+H$ words lie outside triangles and outside private pairs.
Indeed, they have $c_1=1$, a word in a triangle has $c_i=1$ only in its local
columns of profile $(1,1)$, and a word of a private pair has $c_i=1$ only in
its private column. Write $z=n_0$ and $b=n_2$, so that
\begin{equation}\label{eq:partition2}
   n=2m+3t+2e+H+b+z,\qquad E=\Delta-m-t\ge0 .
\end{equation}
Let $\cP$ be the packing of Lemma~\ref{lem:volpairs}: contracted triangles,
enlarged private pairs, and the remaining cylinders. By~\eqref{eq:volpairs}
the volume not covered by $\cP$ is at most $E$.

\subsubsection*{Entries equal to two}
Let $\rho$ be $3b$ minus the total number of entries $c_i=2$ in the $b$
words with $c_1\ge2$. Each of the $2e+H$ words has $c_1=1$ and at most one
entry $2$, since $1+2+2<6$ violates~\eqref{eq:prefix}. Let $\sigma$ be the
number of those without an entry $2$. Let $J$ be the number of entries $2$
in these two groups that come from the larger side of a column of profile
$(a,2)$ with $a\ge3$.

\begin{lemma}\label{lem:twos}
\begin{equation}\label{eq:costsmall}
  b+6z+e+\tfrac53H+\tfrac23\rho+\tfrac45\sigma+\tfrac13J+12r
   +12\bigl(\tfrac ng-1\bigr)\ \le\ 2m .
\end{equation}
\end{lemma}

\begin{proof}
We sharpen both sides of~\eqref{eq:cost}. A word with $c_1\ge2$ has at most
three entries $2$, since $2+2+2+2<12$. If it has two, one or no entries $2$,
then by~\eqref{eq:prefix} the first four terms of its comparison sequence can
be replaced by $(2,2,3,5)$, $(2,3,3,4)$ or $(3,3,3,3)$. This lowers the bound
for $F(v)$ by $1/15$, $1/8$ or $1/6$, which is at least $1/18$ for each entry
$2$ missing to three. For a word with $c_1=1$ and no entry $2$, the start
$(1,3,3,5)$ lowers the bound by $1/15$.

On the column side, a mixed column with sides $a\ge d$ contributes
$1+(a-d)^2/(2ad)$. The excess is $1/4$ for $(2,1)$, at least $2a/9$ for
$(a,1)$ with $a\ge3$, and at least $a/36$ for $(a,2)$ with $a\ge3$. Hence
\[
 h+\tfrac12z+\tfrac14e+\tfrac29H+\tfrac1{36}J\ \le\ \sum_vF(v)\ \le\
 n\Bigl(1-\frac1g\Bigr)+t+\frac{2m+2e+H-b}{12}-\frac{\rho}{18}
 -\frac{\sigma}{15}.
\]
Substituting $h\ge n+t-1+r$ from Lemma~\ref{lem:rank} and multiplying by $12$
gives~\eqref{eq:costsmall}.
\end{proof}

\subsubsection*{Compatible columns}
A mixed column other than a private column is \emph{compatible} if each of
its sides contains either both ends or no end of every private pair. Consider
the columns of the difference matrix (rows indexed by private pairs, entries of
$B_u-B_v$) that correspond to sides containing at most two words after the
contraction. Let $r_2\le r$ be the rank of this submatrix. Let $e_0$ be the
number of compatible $(2,1)$ columns, and $c$ the number of compatible
$(2,2)$ columns containing a word outside private pairs. Columns of profile
$(2,1)$ or $(2,2)$ contain no word of a triangle, so contraction does not
affect them.

\begin{lemma}\label{lem:incompatible}
At most $2r_2$ columns of profile $(2,1)$ and at most $4r_2$ columns of
profile $(2,2)$ are incompatible. Moreover
\begin{equation}\label{eq:twos}
   3b+e+H-\rho-\sigma\ \le\ J+12r_2+4c .
\end{equation}
\end{lemma}

\begin{proof}
The words on the larger side of a $(2,1)$ column are not in private pairs, so
the column is incompatible exactly when its single word $u$ is an end of a
private pair. The singleton side then gives a signed unit vector in the
difference matrix. Unit vectors at different rows are independent, so at most
$r_2$ pairs arise. An end $u$ of a private pair has an entry $1$ and hence at
most one entry $2$, so it is the singleton of at most one $(2,1)$ column.

A side with at most two words gives a vector with at most two nonzero
entries. Choose a basis of $r_2$ such vectors; their supports cover at most
$2r_2$ pairs, and every small side is supported there. Only these pairs can
be separated by a side of a $(2,2)$ column. Each end of such a pair lies in at
most one $(2,2)$ column, again because it has at most one entry $2$. This
gives at most $4r_2$ incompatible $(2,2)$ columns.

The two groups of words have $3b-\rho+(2e+H-\sigma)$ entries $2$ in total. At
most $e$ of them come from $(2,1)$ columns (one per column), and $J$ of them
from $(a,2)$ columns with $a\ge3$. An incompatible $(2,2)$ column contains an
end of a private pair, so it carries at most three of these entries, at most
$12r_2$ altogether. A compatible $(2,2)$ column carries at most four, and
none if all its words lie in private pairs.
\end{proof}

\begin{lemma}\label{lem:uncovered}
$e_0\le E$ and $3e_0+2c\le10E$. Consequently
\begin{equation}\label{eq:uncovered}
   e\le E+2r_2,\qquad c+\tfrac32e\le5E+3r_2 .
\end{equation}
\end{lemma}

\begin{proof}
The consequences follow from $e_0\ge e-2r_2$ (Lemma~\ref{lem:incompatible}).

\emph{First inequality.} Compatible $(2,1)$ columns contain no end of a private
pair. Build a directed graph on the words occurring in them: from each of the
two words on the larger side of such a column draw an arc to its single word.
A word has at most one entry $1$, so out-degrees are at most one, and every
component is either an in-tree or contains exactly one directed cycle. A
$2$-cycle would mean two conflicts between its words. Delete all $e_0$
columns from all words. Cylinders of words not joined by an arc stay disjoint.
The enlarged cylinders of a component remain disjoint from the rest of $\cP$,
because conflicts with outside words do not lie in the deleted columns.
It suffices to show that the union of the enlarged cylinders of a component
built from $a$ columns has volume at least that of the original cylinders
plus $a$.

Every vertex with an outgoing arc has an entry $1$, hence at most one entry
$2$, so it is the singleton of at most one of the columns.

\smallskip\noindent\emph{Trees.} A tree component has $2a+1$ vertices. Let
$d_0$ be the number of its columns whose singleton is the root. After
deletion the root has volume $2^{d_0}$, $a-d_0$ vertices have volume $4$
(singletons of one column) and $a+d_0$ vertices have volume $2$. The sum is
$6a-2d_0+2^{d_0}\ge6a$. Each colour class of the bipartite tree consists of
pairwise disjoint cylinders, and volumes are integers. If one class has sum
greater than $3a$, the union has volume at least $3a+1=(2a+1)+a$. Otherwise
both classes have sum exactly $3a$. If only the root has children, then
$a=d_0$ and the $2a$ leaves alone give $4a\ge3a+1$. If not, take a non-root
internal vertex farthest from the root. It is the singleton of one column, so
it has exactly two children, both leaves. Their enlarged cylinders have volume
$2$, are disjoint, and have distinct supports, because both lost only the
same position. If both lay inside the parent's cylinder of volume $4$, they
would be its two halves, which fix the same extra position, a contradiction.
So one leaf has at least one unit outside its parent. This unit is outside
every other cylinder of the parent's class, since a leaf is adjacent only to
its parent. Hence the union has volume at least $3a+1$.

\smallskip\noindent\emph{Cycles.} A component with a cycle has $2a$ vertices,
each with one outgoing arc, and exactly $a$ of them are singletons of a
column. After deletion the volumes sum to $4a+2a=6a$. Cylinders with a common
point form a clique in the underlying graph. In a graph with out-degrees at
most one, every clique with three vertices is a directed $3$-cycle. If the
cycle has length at least $4$, no point is covered more than twice, and the
union has volume at least $3a=2a+a$.

Let the cycle be $x_1\to x_2\to x_3\to x_1$. The projections of $x_1,x_2,x_3$
have weight $k-2$. They are not all equal: otherwise $x_1,x_2,x_3$ would form
a triangle whose local columns have profile $(2,1)$, contrary to
Lemma~\ref{lem:isolation}. If the three cylinders have no common point, the
covering multiplicity is at most $2$ and we are done. Otherwise their
patterns are pairwise compatible, not all equal, and the triple intersection
has volume $U_3\le2$. Writing $U_j$ for the volume covered exactly $j$ times,
we have $6a=U_1+2U_2+3U_3$, so the union has volume $3a+(U_1-U_3)/2$, and it
suffices to find two units covered exactly once.

Each $x_i$ is the singleton of one column, whose larger side consists of its
cycle predecessor and one vertex outside the cycle; call it the child of
$x_i$. If the child is not a leaf, the deepest internal vertex in its subtree
has two leaf children, and as in the tree case one of them has a unit
outside its parent. That unit is covered exactly once. If the child of $x_i$
is itself a leaf $L_i$, it has such a unit unless $L_i\subseteq x_i$ (we
identify words with their enlarged cylinders). Suppose that $L_i\subseteq
x_i$ and $L_j\subseteq x_j$ for some $i\ne j$. Then $L_i$ fixes exactly one
position $s_i$ beyond $x_i$. As $L_i$ is not adjacent to $x_j$, and the
patterns of $x_i,x_j$ are compatible, $L_i$ conflicts with $x_j$ at $s_i$.
Similarly $x_i$ conflicts with $L_j$ at $s_j$. Here $s_i\ne s_j$ because $x_i$
fixes $s_j$ but not $s_i$. So $L_i$ and $L_j$ conflict at both $s_i$ and $s_j$,
which is impossible. Hence at least two of the three branches contribute a
unit covered exactly once.

\emph{Second inequality.} Return to $\cP$ and treat each of the $e_0+c$
columns separately. For a $(2,1)$ column $\kappa$ with larger side $u,u'$ and
singleton $y$, flip the bit at $\kappa$ in $u$ and $u'$. The flipped cylinders
are disjoint from each other and from all members of $\cP$ except $[y]$. Mark
the set they cover outside $[y]$; its volume is at least $2-1=1$.

For a compatible $(2,2)$ column $\kappa$ without ends of private pairs, let
$U_0,U_1$ be the unions of the cylinders of its $0$-side and of its $1$-side
after deleting $\kappa$. The four words still have distinct supports, so
over $\mathbb F_2$ the indicator of $U_0\triangle U_1$ is a sum of four
multilinear polynomials of degree $k-1$ with distinct leading monomials. It is
therefore a nonzero polynomial of degree $k-1$. A nonzero multilinear
polynomial of degree $j$ is nonzero on a set of measure at least $2^{-j}$: fix
the variables outside one of its monomials of top degree, and sum over the
remaining subcube to get $1$. So $\mu(U_0\triangle U_1)\ge2^{-(k-1)}$. Mark
$\bigl((U_1\setminus U_0)\times\{x_\kappa=0\}\bigr)\cup
\bigl((U_0\setminus U_1)\times\{x_\kappa=1\}\bigr)$, of volume at least $1$.
If one side of a compatible $(2,2)$ column is a private pair, its enlarged
union has volume at least $3$, and the other side consists of two single
words. Flip the pair's union at $\kappa$ and mark the part outside the two
single cylinders, of volume at least $1$. In all cases the members of $\cP$
outside the column have a star at $\kappa$ and keep their conflicts with the
column's words, so the marked sets lie in the region not covered by $\cP$.

Fix a point marked by $h_1$ columns of profile $(2,1)$ and $h_2$ of profile
$(2,2)$. For each of these columns choose a word whose moved cylinder contains
the point. This word (with its private bit deleted, if any) disagrees with
the point exactly at that column. So distinct columns give distinct words.
No two of them are ends of one private pair, because on compatible columns
both ends have equal letters. Two chosen words therefore keep their unique
conflict, and it lies at one of the two columns where they disagree with the
point. At the column of a chosen word, the words with the bit of the point
form the opposite side, which has one word for $(2,1)$ and two words for
$(2,2)$. Writing $\varsigma=h_1+h_2$, we get
$\binom\varsigma2\le h_1+2h_2$. Hence $\varsigma\le5$; if $\varsigma=5$
then $h_1=0$, and if $\varsigma=4$ then $h_1\le2$. In all cases
$3h_1+2h_2\le10$. Integrating over the uncovered region gives
$3e_0+2c\le10E$.
\end{proof}

\begin{proof}[Proof of Theorem~\ref{thm:B}]
We first show $r_2\le r-1$. Dropping nonnegative terms
from~\eqref{eq:costsmall} and using~\eqref{eq:partition2} gives
$n\le4m+3t+e-12r-12(n/g-1)$. Then $e\le E+2r_2\le3E+2r$ and $E=\Delta-m-t$
give
\[
   n\le3\Delta+m-10r-12(n/g-1).
\]
Since $n>3\Delta$, we get $m>10r\ge0$, so $r\ge1$ by Lemma~\ref{lem:rank}. If
$r_2=r$, the column space of the difference matrix would be spanned by $r$
vectors supported on at most $2r$ rows. Then all its rows outside these
would vanish, contradicting that all $m>2r$ rows are nonzero.

Now add $\tfrac{13}2$ times~\eqref{eq:costsmall}, $\tfrac32$
times~\eqref{eq:twos} and $6$ times the second inequality
of~\eqref{eq:uncovered}. The variable $c$ cancels. Substituting
$n=2m+3t+2e+H+b+z$ and $\Delta=m+t+E$, a direct computation gives
\begin{equation}\label{eq:final}
\begin{aligned}
 11n+78\bigl(\tfrac ng-1\bigr)+78r-46r_2+2t+28z+\tfrac43H
   +\tfrac{17}6\rho+\tfrac{37}{10}\sigma+\tfrac23J&\\
   {}+5\,(E+2r_2-e)&\ \le\ 35\Delta .
\end{aligned}
\end{equation}
All terms in the first line except the first four are nonnegative, and so is
$E+2r_2-e$ by~\eqref{eq:uncovered}. Since $r_2\le r-1$ and $r\ge1$, we have
$78r-46r_2\ge32r+46\ge78$. Hence $11n+78n/g\le35\Delta=35(Q-n)$, that is,
$(46+78/g)n\le35Q$. This proves the upper bound for $k\ge4$ when $M_k>g_k$.
If $M_k=g_k$, it holds trivially. The cases $k=2,3$ are covered by the next
subsection, and $M_k\ge g_k$ is given by $G_k$.
\end{proof}

\subsection{Small weights}\label{subsec:smallweights}
For $k=7$ we have $g_7=96$ and $35\cdot128/(46+78/96)<96$, so $M_7=96$. Then
$M_j\le M_7/2^{7-j}=g_j$ for $2\le j\le7$, by~\eqref{eq:Mdouble}. Hence
$M_j=g_j$ for $2\le j\le7$, and the bound of Theorem~\ref{thm:B} also holds
for $k=2,3$.

For $k=8$ the bound leaves the values $192$ and $193$. Suppose that
$n=M_8=193$; then $g=192$ and $\Delta=63$. Subtracting $11n+78(n/g-1)+78$
from $35\Delta$ leaves $115/32$. If $r\ge2$, then $78r-46r_2\ge32r+46$
exceeds $78$ by at least $32$; so $r=1$ and $r_2=0$. In~\eqref{eq:final},
$28z$ and $\frac{37}{10}\sigma$ force $z=\sigma=0$. The integer $E-e\ge0$
appears with coefficient $5$, so $E=e$. What remains is
\begin{equation}\label{eq:slack}
   2t+\tfrac43H+\tfrac{17}6\rho+\tfrac23J\le\tfrac{115}{32}.
\end{equation}
Now $m=63-t-e$ and $b=67-t-H$. Substituting into~\eqref{eq:costsmall} gives
$3e\le47-t-(2H+2\rho+J)/3-1/16$, so $e\le15$. The second inequality
of~\eqref{eq:uncovered} gives $4c\le14e$. With~\eqref{eq:twos} this yields
$201-3t-2H-\rho-J\le13e\le195$, so $3t+2H+\rho+J\ge6$. But~\eqref{eq:slack}
multiplied by $3/2$ gives $3t+2H+\rho+J\le345/64<6$. Hence $M_8=192$, which
completes the proof of Theorem~\ref{thm:B}.

\begin{proof}[Proof of Theorem~\ref{thm:C}]
By Lemma~\ref{lem:coding} and Theorem~\ref{thm:B},
$S_d\le M_{d+1}<\frac{35}{46}\,2^{d+1}$, and $S_d\le M_{d+1}=g_{d+1}$ for
$d\le7$. The lower bounds come from Theorem~\ref{thm:A}.
\end{proof}

Theorem~\ref{thm:B} gives $384\le M_9\le387$ and $768\le M_{10}\le777$. Since
$M_k/2^k$ is nondecreasing, its limit $\lambda$ exists and satisfies
$\frac34\le\lambda\le\frac{35}{46}=\frac34+\frac1{92}$.

\section{Open problems}\label{sec:open}
\begin{enumerate}[1.]
\item Determine $S_4\in[18,24]$. Is the extremal $5$-code $G_5$ unique, so
that $S_4\le23$ (Proposition~\ref{prop:S4})?
\item What is the largest size $b_d$ of a based, pierced, neighbourly family
of $d$-simplices in $\R^d$? We know $b_3\le4$ (Corollary~\ref{cor:based3}),
$9\le b_4\le M_4=12$
(Remark~\ref{rem:basedcodes}), and $S_d\ge2b_d$. Does $b_{d+1}/b_d$ exceed
$2$ for larger $d$? Determine $\lim_d S_d/2^d\in[\frac98,\frac{35}{23}]$, if
the limit exists.
\item Prove Conjecture~\ref{conj:34}, that is, $M_k=3\cdot2^{k-2}$; it would
give $S_d\le\frac32\,2^d$. For $V$ a neighbourly $k$-code, let $\tau(V)$ be
the least number of words meeting all rainbow triples. By
Lemma~\ref{lem:split}, the inequality $\tau(V)\le\Delta(V)$ for one code of
size $M_k$ in every weight would imply the conjecture. It holds for $G_k$,
with equality, and for codes with a linear rainbow hypergraph
(Theorem~\ref{thm:linear}). Does it hold for all neighbourly codes?
\item Is $G_k$ the unique neighbourly $k$-code of size $g_k$ for
$6\le k\le8$?
\end{enumerate}

\appendix
\section{Vertices of the nine simplices}\label{app:vertices}
Table~\ref{tab:vertices} lists the vertices of $P_1,\dots,P_9$. The vertex
$w_0$ is opposite to the base $\{t=0\}$, and $w_j$ ($j=1,\dots,4$) is opposite
to the $j$-th side facet in the order of Table~\ref{tab:words}. Vertices of
the $18$ simplices of Remark~\ref{rem:explicit} are obtained by applying the
maps $A$ and $B$.

\begin{table}[ht]
\centering\footnotesize
\setlength{\tabcolsep}{5pt}
\renewcommand{\arraystretch}{1.05}
\caption{Vertices of $P_1,\dots,P_9$.}\label{tab:vertices}
\begin{tabular}{@{}ccrrrr@{}}
\toprule
 & & $x$ & $y$ & $z$ & $t$\\
\midrule
$P_{1}$ & $w_0$ & $368847/277$ & $1002/277$ & $466/831$ & $12302/831$ \\
 & $w_1$ & $-4762/1659$ & $50/1659$ & $-97/553$ & $0$ \\
 & $w_2$ & $-31723/1461$ & $74/487$ & $-1046/1461$ & $0$ \\
 & $w_3$ & $4119/881$ & $-144/881$ & $346/881$ & $0$ \\
 & $w_4$ & $2889/488$ & $11/488$ & $59/488$ & $0$ \\
\addlinespace[2pt]
$P_{2}$ & $w_0$ & $-22267/33$ & $-2020/99$ & $2654/165$ & $22906/495$ \\
 & $w_1$ & $80111/1610$ & $583/1610$ & $-48/161$ & $0$ \\
 & $w_2$ & $5191/53$ & $830/583$ & $-1312/583$ & $0$ \\
 & $w_3$ & $6469/119$ & $-113/238$ & $61/34$ & $0$ \\
 & $w_4$ & $20040/551$ & $-21/551$ & $-22/551$ & $0$ \\
\addlinespace[2pt]
$P_{3}$ & $w_0$ & $3475449/88378$ & $13943/176756$ & $-6893/353512$ & $181715/353512$ \\
 & $w_1$ & $577072/16069$ & $917/16069$ & $4224/16069$ & $0$ \\
 & $w_2$ & $428906/10451$ & $1777/10451$ & $576/10451$ & $0$ \\
 & $w_3$ & $839054/35813$ & $2824/35813$ & $1713/35813$ & $0$ \\
 & $w_4$ & $23398/1651$ & $-2077/1651$ & $-1839/1651$ & $0$ \\
\addlinespace[2pt]
$P_{4}$ & $w_0$ & $2900402/76741$ & $2935/76741$ & $1826/76741$ & $28353/76741$ \\
 & $w_1$ & $-350/779$ & $-1669/779$ & $-1589/779$ & $0$ \\
 & $w_2$ & $185960/9781$ & $-6761/9781$ & $-661/9781$ & $0$ \\
 & $w_3$ & $-1382/35$ & $13/7$ & $1/5$ & $0$ \\
 & $w_4$ & $75862/2567$ & $-841/2567$ & $-359/2567$ & $0$ \\
\addlinespace[2pt]
$P_{5}$ & $w_0$ & $137462008/3449981$ & $367398/3449981$ & $21407/6899962$ & $2472415/6899962$ \\
 & $w_1$ & $428906/10451$ & $1777/10451$ & $576/10451$ & $0$ \\
 & $w_2$ & $3147/64$ & $8023/22976$ & $-197/718$ & $0$ \\
 & $w_3$ & $-17459/8287$ & $-3537/66296$ & $2471/66296$ & $0$ \\
 & $w_4$ & $1317986/37017$ & $-4399/37017$ & $-6713/37017$ & $0$ \\
\addlinespace[2pt]
$P_{6}$ & $w_0$ & $1812959/54426$ & $-6539/54426$ & $2425/54426$ & $17749/36284$ \\
 & $w_1$ & $-38381/45$ & $-1198/45$ & $154/9$ & $0$ \\
 & $w_2$ & $203507/1644$ & $4747/1644$ & $-1451/1644$ & $0$ \\
 & $w_3$ & $252781/210$ & $3121/105$ & $-2323/210$ & $0$ \\
 & $w_4$ & $14118/463$ & $-99/463$ & $34/463$ & $0$ \\
\addlinespace[2pt]
$P_{7}$ & $w_0$ & $18641377/333942$ & $31376/166971$ & $-15921/111314$ & $5177/7951$ \\
 & $w_1$ & $-4762/1659$ & $50/1659$ & $-97/553$ & $0$ \\
 & $w_2$ & $379759/37212$ & $-506/9303$ & $2477/12404$ & $0$ \\
 & $w_3$ & $58787/20454$ & $-1199/10227$ & $1753/6818$ & $0$ \\
 & $w_4$ & $56531/6720$ & $137/6720$ & $459/2240$ & $0$ \\
\addlinespace[2pt]
$P_{8}$ & $w_0$ & $50159269/1501968$ & $4305/62582$ & $38245/1501968$ & $124291/375492$ \\
 & $w_1$ & $147431/43608$ & $-56/5451$ & $57/14536$ & $0$ \\
 & $w_2$ & $379759/37212$ & $-506/9303$ & $2477/12404$ & $0$ \\
 & $w_3$ & $1106239/7608$ & $2101/1902$ & $-2179/2536$ & $0$ \\
 & $w_4$ & $-50975/404$ & $569/101$ & $241/404$ & $0$ \\
\addlinespace[2pt]
$P_{9}$ & $w_0$ & $3295442/92375$ & $-3839/92375$ & $1878/92375$ & $83094/92375$ \\
 & $w_1$ & $80111/1610$ & $583/1610$ & $-48/161$ & $0$ \\
 & $w_2$ & $779464/17713$ & $4169/17713$ & $-1152/17713$ & $0$ \\
 & $w_3$ & $490213/9772$ & $704/2443$ & $-1107/9772$ & $0$ \\
 & $w_4$ & $14118/463$ & $-99/463$ & $34/463$ & $0$ \\
\bottomrule
\end{tabular}
\end{table}
\clearpage

\section*{Declaration of generative AI and AI-assisted technologies}
During the preparation of this work the authors used ChatGPT Astra Pro~6.2 
(OpenAI) and Claude (Anthropic, model Claude Opus~5.5). ChatGPT Astra found
the construction of Theorem~\ref{thm:A}: the nine simplices of
Section~\ref{sec:nine} and the explicit family of $18$ simplices of
Remark~\ref{rem:explicit}. Claude independently verified the
construction in exact rational arithmetic, including the family of $36$
simplices in $\R^5$ of Remark~\ref{rem:explicit}, checked the proofs, merged
the authors' drafts into a single manuscript, edited the text, and
contributed Corollary~\ref{cor:based3}, Remarks~\ref{rem:heights}
and~\ref{rem:basedcodes}, and Example~\ref{ex:rainbow18}(b). The authors
reviewed and verified all of the content and take full responsibility for
this publication.


\end{document}